\documentclass[11pt]{article}
\usepackage{amsmath,amssymb}
\newcommand{\sect}[1]{\setcounter{equation}{0}\section{#1}}
\newcommand{\subsect}[1]{\subsection{#1}}

\newtheorem{lemma}{Lemma}[section]

\newfont{\extra}{msbm10 scaled\magstep1}
\newcommand{\extr}[1]{\mbox{\extra #1}}
\def\C{\extr C}

\begin{document}

\begin{center}
{\LARGE {\bf {QUANTIZATION  \\[1mm]
OF LIE-POISSON STRUCTURES \\[3mm]
BY  PERIPHERIC CHAINS}}} \vskip1cm

Vladimir D. Lyakhovsky $^1$  and Mariano A. del Olmo $^2$
 \vskip0.5cm

$^1$ {\it Theoretical Department, St. Petersburg State University, \\[0.15cm]
198904, St. Petersburg, Russia \\[0.25cm]
$^2$ Departamento de F\'{\i}sica Te\'orica, Universidad de
Valladolid, \\[0.15cm]
E--47011, Valladolid,  Spain}\\[0.3cm]

E-mail: lyakhovs@pobox.spbu.ru , olmo@fta.uva.es
\end{center}

\vskip1cm

\centerline{\today}

\vskip1.5cm

\begin{abstract}
The quantization properties of composite peripheric twists are studied.
Peripheric chains of extended twists are constructed for $U(sl(N))$ in
order to obtain composite twists with
sufficiently large carrier subalgebras. It is proved that the peripheric
chains can be enlarged with  additional Reshetikhin and Jordanian factors.
This provides the
possibility to construct new solutions to Drinfeld equations and, thus,
to quantize new sets of Lie-Poisson
structures. When the Jordanian additional factors are used the carrier
algebras of the enlarged peripheric chains are transformed into   algebras of
motion of the form
${\mathfrak g}_{{\cal JB}}^{{\cal P}}={\mathfrak g}_{{\cal H}}\vdash
{\mathfrak g}_{{\cal P}}$.  The factor algebra ${\mathfrak g}_{{\cal H}}$
is a direct sum of
Borel and contracted Borel  subalgebras of lower dimensions. The
corresponding $\omega$--form is a coboundary. The enlarged peripheric chains
${\cal F}_{{\cal JB}}^{{\cal P}}$ represent the twists that
contain operators external with respect to the Lie-Poisson structure.
The properties of new twists are illustrated by quantizing
$r$--matrices for  the
algebras $U(sl(3))$,  $U(sl(4))$ and $U(sl(7))$.
\end{abstract}

\newpage
\sect{Introduction}
\label{Introduction}

Quantum algebras are often considered
as deformations of symmetries described by Lie algebraic structures.
Many of their physical applications are dealing with composite systems
\cite{cesargomez}--\cite{pasquier}. In these cases the coproducts
(i.e, the coalgebraic  structure called
the co-structure) are
relevant \cite{angela}.
When composite systems (i.e. the sets composed of sub-systems) are studied
physicists are well accustomed to work with Lie algebras
and with nondeformed Hopf algebras (symmetric or cocommutative coalgebras)
since the set of sub-systems is not ordered and for that  all
measurable results must be independent of the order  chosen for  the
sub-systems (remember the case of the composition of angular momenta).

In general, in a quantum algebra both structures, algebraic and coalgebraic,
are deformed. In other words, the coproduct in non cocommutative, thus,
quantum algebras are ordered structures. We can interpret this lack of
symmetry as an existence of a kind of correlation between the sub-systems.
In this way, quantum groups supply us with an approximate symmetry
\cite{olmo02}.

Twisted quantum algebras have the nice property:
one can always conserve the classical Lie algebraic compositions (use the
nondeformed basis), only the co-structure is deformed. This is why twisted
quantum algebras are interesting objects when physical applications are
considered.

It was shown by Drinfeld  \cite{D83}
that a quasitriangular Hopf algebra $\cal A$ can be
transformed into a twisted
Hopf algebra when a certain adjoint operator of the
``twisting'' element of $ {\cal F} \in
{\cal A}\otimes {\cal A}$ is applied to the coproducts of  $\cal A$.
The algebraic structure of the twisted algebra
 ${\cal A}_{\cal F}$ is the same as in
the original one but the co-structure in ${\cal A}_{\cal F}$  and
the quantum $\cal R$--matrix are different.
Thus, using twists we can obtain quantum deformations of universal
enveloping algebras of Lie algebras where the irreducible representations are
classical but their compositions (tensor products) are deformed.

Since the seminal paper of Drinfeld mentioned above and especially
during the last years important studies were performed in the
area of twist deformations and complicated twisting elements defined on large
``carrier" algebras were constructed \cite{RES}--\cite{aizawa}.
Notwithstanding their complicated form
one can always decompose them into the product of
basic twisting elements (factors) \cite{LYA}: Reshetikhin  \cite{RES},
Jordanian
\cite{OGI} or extended \cite{KLM} factors. The twists are often constructed
in a form of multiparametric families and the
boundary twisting elements have the specific properties. In this concern
we can mention the peripheric extended twists which are the limit
versions of the extended jordanian twists \cite{LO}.

In this work we study the possibility
to construct chains of peripheric extended twists. We
shall prove that these chains can be enlarged using
the Reshetikhin-like or the Jordanian-like factors.
The enlarged chains of twists provide the possibility  to deform universal
enveloping Lie algebras in a way different from the canonical twisting. The
overall adjoint transformation cannot be reduced to the adjoint operators of the
carrier algebra of the classical $r$-matrix. In other words, from a more physical
point of view, we can quantize new classical $r$--matrices, or equivalently, new
classical systems described by Poisson--Lie structures \cite{charipressley}.

The paper is organized as follows. Section~\ref{preliminaries}
presents a review of the
main facts related
to the theory of twists that we shall use along the paper. In the next
section, using the particular case of the Lie algebra $sl(3)$, we demonstrate
that the carrier algebra (i.e. the minimal subalgebra
of the Lie algebra  where the twisting element ${\cal F}$ can be defined)
associated to a peripheric twist enlarged by a Reshetikhin or Jordanian
factor is different from the initial carrier algebra of the peripheric twist.
In Section~\ref{periphericchains} we prove the existence of
peripheric chains for
$sl(N)$, and also find the maximal number of peripheric twisting compositions
(or links) that we can add to
construct a peripheric chain for $sl(N)$ in the case of even $N=2n$ or odd
$N=2n-1$ dimension.
The properties of peripheric chains for $U(sl(N))$ are studied in the
following section.
We discuss separately two cases depending on the character of the additional
twisting factors: Reshetikhin  or
Jordanian. In the first case the enlarged peripheric chains allows to
quantize new
classical $r$--matrices whose associated $\omega$--forms are cohomologically
nontrivial.
In the second one, it is proved that the carrier
algebras of the peripheric chains enlarged by the Jordanian factors have
the structure of semidirect sum of Borel and Abelian subalgebras
that can be seen
as motion algebras. Section~\ref{examples}
presents the enlarged peripheric twists for $U(sl(3))$, $U(sl(4))$
and $U(sl(7))$ in order to
illustrate the theory developed in the preceding sections. The discussion
of the obtained results and some comments conclude the paper.

\sect{Mathematical preliminaries}
\label{preliminaries}

Quasitriangular Hopf algebras ${\cal A}(m,{\mit\Delta} ,\eta ,\epsilon ,S,
{\cal R})$ can be transformed \cite{D83} by an invertible twisting
element ${\cal F}\in {\cal A}\otimes {\cal A}$ satisfying the equations
\begin{equation}  \label{drinf}
{\cal F}_{12}({\mit\Delta} \otimes id)({\cal F})={\cal F}_{23}
(id\otimes  {\mit\Delta} )({\cal F}),\quad (\epsilon \otimes id)({\cal F})=
(id\otimes \epsilon)({\cal F})=1 ,
\end{equation}
into twisted Hopf algebras ${\cal A}_{{\cal F}}
( m,{\mit\Delta}_{ {\cal F}},\epsilon ,S_{{\cal F}}, {\cal R}_{{\cal F} })
$, that has the same multiplication and counit but the twisted coproduct,
antipode and
${\cal R}$--matrix defined as follows
\begin{equation}\label{def-t}
\begin{array}{rclcrllr}
{\mit\Delta} _{{\cal F}}(X) & = & {\cal F}{\mit\Delta} (X){\cal F}^{-1},
\qquad  &S_{{\cal F}}(X)\; =\; vS(X)v^{-1},
&\qquad   {\cal R}_{{\cal F}}&=&
\left( {\cal F}\right)_{21}{\cal RF}^{-1}, \\[2mm]
{\cal F} & = & \sum f_{i}^{(1)}\otimes f_{i}^{(2)},
\qquad & X\in {\cal A},
&\qquad  v&=&\sum f_{i}^{(1)}S(f_{i}^{(2)}).
\end{array}
\end{equation}

Twists are mostly used as tools to construct quantum deformations,
${\cal F}:U({\mathfrak g})\to U_{{\cal F}}({\mathfrak g})$, of universal
enveloping algebras $U({\mathfrak g})$ of simple Lie algebras
${\mathfrak g}$,
considered as Hopf algebras with primitive generators $L\in{\mathfrak g}$
(i.e. ${\mit\Delta} (L)=L\otimes 1+1\otimes L$). The minimal subalgebra
${\mathfrak g}_{c}\subset {\mathfrak g}$ on which the twist ${\cal F}$
can be defined is called the carrier of the twist ${\cal F}$.

In general, the composition of twists is not a twist. But there are some
important examples of the opposite behaviour. In particular, a composition of
twists ${\cal F}_{1}$ and ${\cal F}_{2}$ (with carriers
${\mathfrak g}_{c_1}$ and
${\mathfrak g}_{c_2}$) can be performed when in the twisted algebra
$U_{{\cal F}_{1}}( g) $ one can find the primitive generators for
the carrier ${\mathfrak g}_{c_2}$.

There are three basic twisting factors (BTF's) \cite{LYA}. Up to now all the
twists that are known in the explicit form are composed of these BTF's:
\medskip

{\bf 1.-} The Reshetikhin basic twisting factor (RF) is defined for any pair
of commuting primitive elements $l_{i},l_{j}$ in ${\cal A}$ \cite{RES} by
\begin{equation}
{\mit\Phi} _{R}=e^{\psi _{i,j}l_{i}\otimes l_{j}}.  \label{resh}
\end{equation}
The coefficient $\psi _{i,j}$ plays the role of the deformation parameter.
\medskip

{\bf 2.-}  The Jordanian BTF (JF) \cite{OGI} has as  a carrier the
two-dimensional (2D) Borel algebra ${\bf B}$
\begin{equation}
\lbrack H,E]=E.  \label{bor}
\end{equation}
For $U({\bf B})$ the following element is the solution of (\ref{drinf}):
\begin{equation}
{\mit\Phi} _{J}=e^{H\otimes \sigma }  \label{jord}
\end{equation}
with
\begin{equation}
\sigma =\ln (1+E).  \label{sigm}
\end{equation}
The result of the deformation ${\mit\Phi} _{J}:U({\bf B})\to U_{J}(
{\bf B})$ is a Hopf algebra $U_{J}({\bf B})$ with the coproducts:
\begin{equation}\label{jord-co}
\begin{array}{l}
\Delta _{J}(H)=H\otimes e^{-\sigma }+1\otimes H, \\[0.2cm]
\Delta _{J}(E)=E\otimes e^{\sigma }+1\otimes E.
\end{array}
\end{equation}
The deformation parameter is introduced by scaling the element $E\in {\bf B}$:
\begin{equation}\label{jord-xi}
\begin{array}{llllll}
& E& \longrightarrow  & \xi E, & &\\[0.2cm]
& \sigma & \longrightarrow &\sigma (\xi ) &=&\ln (1+\xi E), \\[0.2cm]
& {\mit\Phi}_{J} &\longrightarrow  &{\mit\Phi} _{J}(\xi )
&=&e^{H\otimes \sigma (\xi )}.
\end{array}
\end{equation}

{\bf 3.-}  The third BTF \cite{KLM}, the  extending factor (EF) or simply
the  extension, is defined for the deformed universal enveloping
algebra $U_{J}({\bf H})$, where ${\bf H}$ is the 3D
Heisenberg algebra
\begin{equation} \label{heis}
\lbrack A,B]=E,\qquad \lbrack A,E]=0, \qquad [B,E]=0,
\end{equation}
with the coproducts
\begin{equation}\label{heis-co}
\begin{array}{l}
\Delta _{J}(A)=A\otimes e^{\alpha \sigma }+1\otimes A, \\[0.2cm]
\Delta _{J}(B)=B\otimes e^{\beta \sigma }+1\otimes B,\qquad \alpha +\beta =1 ,
\\[0.2cm]
\Delta _{J}(E)=E\otimes e^{\sigma }+1\otimes E,
\end{array}
\end{equation}
and $\sigma $  as in (\ref{sigm}). The extending factors
\begin{equation}
\label{ext}{\mit\Phi} _{E}=e^{A\otimes Be^{-\beta \sigma }}, \qquad
{\mit\Phi} _{E}^{\prime }=e^{-B\otimes Ae^{-\alpha \sigma }}
\end{equation}
are the solutions of the twist equations (\ref{drinf}) and
induce the deformations
${\mit\Phi} _{E}:U_{J}({\bf H})\to U_{E}({\bf H})$ or
${\mit\Phi} _{E}^{\prime}:U_{J}
({\bf H})\to U_{E^{\prime }}({\bf H})$. For example,
the co-structure of the algebra
$U_{E}({\bf H})$ is defined by the coproducts:
\begin{equation}\label{ext-co-ur}
\begin{array}{l}
\Delta _{E}(A)=A\otimes e^{-\beta \sigma }+1\otimes A, \\[0.2cm]
\Delta _{E}(B)=B\otimes e^{\beta \sigma }+e^{\sigma }\otimes B, \\[0.2cm]
\Delta _{E}(E)=E\otimes e^{\sigma }+1\otimes E.
\end{array}
\end{equation}
EF's were first considered as parts of the so-called extended
Jordanian twists (EJ's) \cite{KLM}. Contrary to the other two BTF's these
are the discrete solutions of the twist equations, though they can borrow a
continuous parameter from a smooth curve of equivalent algebras
$U_{J}({\bf B };\sigma (\xi ))$, where $\sigma (\xi )$ is the same
as in the parameterized Jordanian twist (see (\ref{jord-xi})).

Note that usually the EF is applied to the universal
enveloping algebras $U_J({\bf L}(\alpha, \beta))$ with the 4D
carrier subalgebra ${\bf L} (\alpha, \beta)$. The parameters $\alpha$ and
$\beta$ that appear in (\ref{heis-co}) are to coincide with the corresponding
eigenvalues of the operator ${\rm ad}(H)$:
\begin{equation}  \label{carr-4}
\begin{array}{lllll}
&[ H,E ] = E, \qquad &[ H,A ] = \alpha A, \qquad &[ H,B ] = \beta B,& \\[2mm]
&[ A,B ] = E, \qquad &[ E,A ] = 0,\qquad &[ E,B ] = 0,
\qquad &\alpha + \beta = 1.
\end{array}
\end{equation}
This Hopf algebra is already deformed by the Jordanian factor. The two
factors are then considered together forming the extended Jordanian
twist \cite{KLM}, ${\cal F}_{EJ} ={\mit\Phi}_{E}{\mit\Phi}_{J} $. This twist
defines the deformed Hopf algebras $U_{{\cal E}} ( {\bf L}(\alpha,\beta) ) $
with the co-structure
\begin{equation}\label{e-costr}
\begin{array}{lcl}
{\mit\Delta} _{{\cal E}}\,(H) & = & H\otimes e^{-\sigma }+1\otimes H-A\otimes
Be^{-(\beta +1)\sigma }, \\[0.2cm]
{\mit\Delta} _{{\cal E}}\,(A) & = & A\otimes e^{-\beta \sigma }+1\otimes A,
\\[0.2cm]
{\mit\Delta} _{{\cal E}}\,(B) & = & B\otimes e^{\beta \sigma }
+e^{\sigma}\otimes B, \\[0.2cm]
{\mit\Delta} _{{\cal E}}\,(E) & = & E\otimes e^{\sigma }+1\otimes E.
\end{array}
\end{equation}

Also we would like to remind that any BTF can be deformed by the previous
factors of the sequence of twists. Nevertheless, it always
has the form of an exponent whose argument is a tensor product of two
elements of ${\cal A}$. Moreover, when in the sequence of twists all the
parameters independent from the proper parameter of the BTF are equal to zero
the BTF recovers its canonical form.

In \cite{LO} the extended Jordanian twists were studied in the limit points
$\alpha =0$ (or $\beta =0$). The corresponding twists were called
peripheric extended twists (PET's). It is easy to see from (\ref{ext-co-ur})
that in these limiting points some generators of ${\bf L}^{{\cal  P}} =
 {\bf L} (1,0)$ remain primitive  in the deformed
$U_{{\cal  P}}({\bf H})$. For example, in
$U_{{\cal  P}}({\bf L}^{{\cal  P}})$
we have
\begin{equation}\label{pet-costr}
\begin{array}{lcl}
{\mit\Delta}^{{\cal P}}_{{\cal E}}\,(H) & = & H\otimes e^{-\sigma }
+1\otimes H-A\otimes Be^{- \sigma }, \\[0.2cm]
{\mit\Delta}^{{\cal P}}_{{\cal E}}\,(A) & = & A\otimes 1+1\otimes A,
 \\[0.2cm]
{\mit\Delta}^{{\cal P}}_{{\cal E}}\,(B) & = & B\otimes 1
+e^{\sigma}\otimes B,\\[0.2cm]
{\mit\Delta}^{{\cal P}}_{{\cal E}}\,(E) & = &
E\otimes e^{\sigma }+1\otimes E.
\end{array}
\end{equation}
The second cohomology group $H^{2}\left( {\bf L}^{{\cal P}}\right) $ of the
carrier algebra for the PET is nontrivial.
We have $\dim H^{2}\left( {\bf L}^{{\cal P}}\right) =1$ (in contrast to the
case $\alpha ,\beta \neq 0,1$). The nontrivial cocycle can be chosen to be
proportional to $H^{\ast }\wedge A^{\ast }$. This means that
${\bf L}^{{\cal P}}$ has not only the nondegenerate coboundary
$\omega =E^{\ast }\left( \left[ ,\right] \right) $ but\ also the
nondegenerate cocycle \cite{STO}
\begin{equation}\label{form-re}
\omega _{{\cal RE}}=E^{\ast }\left( \left[ ,\right] \right)
+\xi H^{\ast }\wedge A^{\ast }.
\end{equation}
The latter defines the classical $r$--matrix
$r_{{\cal RE}}=H\wedge E+A\wedge B+\psi A\wedge E$. It is easy
to verify that the
corresponding twisting element is a composition of
the PET ${\cal F}_{{\cal E}}^{{\cal P}}$
and the RF
\begin{equation}
\mit\Phi _{R}=e^{\psi A\otimes \sigma };  \label{resh-add}
\end{equation}
i.e.
\begin{equation}\label{tw-pre}
{\cal F}_{{\cal RE}}^{{\cal P}}=\mit\Phi _{R}{\cal F}_{{\cal E}}^{{\cal P} }
=e^{\psi A\otimes \sigma }e^{A\otimes B}e^{H\otimes \sigma }.
\end{equation}
Other specific possibilities of PET's
arise when ${\bf H}(\beta =0)$ is considered as
a subalgebra of a simple Lie algebra $\mathfrak g$.
The most interesting are the
non-Abelian primitive subalgebras in $\mathfrak g$
that contain the generators like $A$
(see Section 4 where these possibilities are studied in detail).
In the simplest case when the initial twist is a peripheric twist
with several extension factors ${\cal F}_{{\cal E}}^{{\cal P}}$
these additional
possibilities lead to the deformations equivalent to the
multi-Jordanian twistings (see Section~\ref{examples}).

More fruitful and much more complicated is the situation where the initial
twist is a nontrivial composition of extended twists.
For Hopf algebras $U({\mathfrak g})$ with classical simple Lie algebras
$\mathfrak{g}$
there exists the possibility to construct systematically the compositions
of twists called chains \cite{KLO}:
\begin{equation}
\label{mainchain}
{\cal F}_{{\cal B}_{p\prec 0}}\equiv {\cal F}_{{\cal B}_{p}}
{\cal F}_{{\cal B}_{p-1}}\ldots {\cal F}_{{\cal B}_{0}}.
\end{equation}
The factors ${\cal F}_{{\cal B}_{k}}={\mit\Phi} _{{\cal E}_{k}}
{\mit\Phi} _{{\cal J}_{k}} $ of the canonical chain are the canonical
extended Jordanian twists for the initial Hopf algebra ${\cal A}$.
Their carriers
are the multidimensional analogs to ${\bf L}( {1}/{2},{1}/{2}) $,
the extensions
$\{ {\mit\Phi} _{{\cal E}_{k}},k=0,\ldots , p-1 \} $ contain the fixed set of
normalized factors like ${\mit\Phi} _{{\cal E }}=
\exp \{A\otimes Be^{-\frac{1}{2}\sigma}\}$.
The sequence of twists in a chain for the Lie algebra ${\cal A}=U(sl(N))$
corresponds to the sequence of injections
$U(sl(N))\supset U(sl(N-2))\supset \ldots \supset U(sl(N-2k))\ldots $.
For each ${\cal A}_k$
the  initial root $\lambda _{0}^{k}$ is fixed. The extensions
${\mit\Phi} _{{\cal E}_{k}}$ are
defined by the set $\pi _{k}$ of constituent roots
$\pi _{k}= \{ \lambda ^{\prime },\;
\lambda ^{\prime \prime }\;|\; \lambda^{\prime }
+\lambda^{\prime\prime }=\lambda _{0}^{k};\; \lambda ^{\prime}
+\lambda _{0}^{k},\;
\lambda^{\prime \prime }+\lambda _{0}^{k}\notin \Lambda_{{\cal A}} \} $,
where $\Lambda _{{\cal A}}$ is the root system of ${\cal A}_{0}$.
The set of
constituent roots can be naturally decomposed as
$\pi _{k}=\pi _{k}^{\prime}\, \cup
\,\pi _{k}^{\prime \prime},\  \pi _{k}^{\prime }=\{\lambda ^{\prime }\},\,
\pi_{k}^{\prime \prime }=\{\lambda ^{\prime \prime }\}$. In these terms
the factors
${\mit\Phi} _{{\cal E}_{k}}$ and
${\mit\Phi}_{{\cal J}_{k}}$ have the form
\begin{equation}
{\mit\Phi} _{{\cal J}_{k}}=\exp \{H_{\lambda _{0}^{k}}\otimes \sigma
_{0}^{k}\},\qquad \quad \sigma _{0}^{k}=
\ln (1+L_{{\lambda_{0}^{k}}});
\end{equation}
\begin{equation}
\label{gen-ext}
{\mit\Phi} _{{\cal E}_{k}}=\prod_{\lambda ^{\prime } \in \pi _{k}^{\prime }}
{\mit\Phi} _{{\cal E}_{\lambda ^{\prime }}}=\prod_{\lambda ^{\prime }
\in \pi_{k}^{\prime}}\exp \{L_{\lambda ^{\prime }}
\otimes L_{
{\lambda _{0}^{k}} -\lambda ^{\prime }}
e^{-\frac{1}{2}\sigma _{0}^{k}}\},
\end{equation}
where the generator $L_{\lambda }$  is associated to the root $\lambda $.

The role of additional twistings provided by the
chains of twists with the peripheric properties and the new class of
Lie-Poisson structures that can be thus quantized are studied in
Section~\ref{propertieschains}.

\sect{Deformations of carriers: the $sl(3)$ case}
\label{deformationcarriers}

Most interesting are the cases where the additional factors change the carrier
subalgebra. This may happen if the initial carrier ${\bf L}^{{\cal P}}$ is a
proper subalgebra  ${\bf L}^{{\cal P}} \subset {\mathfrak g}$ and the Drinfeld
equations are considered over  ${\mathfrak g}$. In this Section the simplest
variant of such situation is studied.
Suppose that  $\mathfrak g$ contains the Lie subalgebra
${\bf M}=H^{\perp }\triangleright {\bf L} ^{{\cal P}}$, which
is the extension
of ${\bf L}^{{\cal P}}$ by the 1D subalgebra generated by  $H^{\perp }$
with the following action on ${\bf L}^{{\cal P}}$:
\[
\lbrack H^{\perp },E]=0,\qquad \lbrack H^{\perp },A]=A,\qquad \lbrack
H^{\perp},B]=-B,\qquad \lbrack H^{\perp },H^{\cal P}]=0.
\]
When ${\mathfrak g}=sl(3)$ the generators of ${\bf M} \subset sl(3) $
can be identified as follows
\begin{equation}\label{generat-3}
\begin{array}{llll}
 & H^{\cal P}=\frac{1}{3}\left( 2E_{11}-E_{22}-E_{33}\right),&
\quad A=E_{12}, &\quad E=E_{13} ,\\[0.3cm]
&H^{\perp }=\frac{1}{3}\left( E_{11}-2E_{22}+E_{33}\right),
&\quad B=E_{23} .&
\end{array}\end{equation}
In $U_{{\cal P}}( {\bf M}) $ twisted by ${\cal F}_{{\cal E}}^{{\cal P}}$
the co-structure is defined by
(\ref{pet-costr}) while the coproduct of the additional
generator $H^{\perp }$ remains primitive,
${\mit\Delta} _{{\cal  P}}( H^{\perp})
=H^{\perp }\otimes 1+1\otimes H^{\perp}$.
If we want to change the carrier
${\bf L}^{{\cal P}}$ of the twist by introducing an extra factor (to enlarge
the ${\cal F}_{{\cal E}}^{{\cal P}}$) we have to consider
the primitive Borel subalgebra ${\bf B} \subset {\bf M}$ generated by
$\{ A,H^{\perp }\} $ or the primitive Abelian subalgebra ${\bf A}$
generated by
$\{\sigma ,H^{\perp }\} $. This means that $U_{{\cal P}}( {\bf M})$
admits additional Jordanian or Reshetikhin twists
\[
{\mit\Phi} _{{\cal J}}^{\perp }=\exp \{H^{\perp }\otimes
\sigma_{A}(\xi )\} \quad
{\rm or} \quad
{\mit\Phi} _{{\cal R}}=\exp \{\zeta H^{\perp}\otimes \sigma \}.
\]
Composing ${\mit\Phi} _{{\cal J}}^{\perp }$ or ${\mit\Phi} _{{\cal R}}$ with
${\cal F}_{{\cal E}}^{{\cal P}}$ we obtain new solutions of the Drinfeld
equations for the algebra $U(sl(3))$:
\begin{equation}
\label{simp-add-jord}
\begin{array}{rcl}
{\cal F}_{{\cal J}{\cal E}}^{{\cal P}} & = & {\mit\Phi} _{{\cal J}
}^{\perp } {\cal F}_{{\cal E}}^{{\cal P}}=e^{H^{\perp }\otimes \sigma_{A}(\xi)}
e^{A\otimes B}e^{H^{{\cal P}} \otimes \sigma }, \\[0.2cm]
{\cal F}_{{\cal J}{\cal E}}^{{\cal P}} & : & U\left( {\bf M}\right)
\to U_{{\cal P}}\left( {\bf M}\right) \to U_{{\cal JP}}
 \left( {\bf M}\right) ;
 \end{array}
\end{equation}
and
 \[
\begin{array}{rcl}
{\cal F}_{{\cal RE}} & = & {\mit\Phi}_{{\cal R}}
{\cal F}_{{\cal E}}^{{\cal P}}=
e^{\zeta H^{\perp }\otimes \sigma }e^{A\otimes B}e^{H^{{\cal P}}
\otimes \sigma },\\[0.2cm]
{\cal F}_{{\cal RE}} & : & U\left( {\bf M}\right) \to U_{{\cal E}}
\left( {\bf M}\right) \to U_{{\cal RE}}\left( {\bf M}\right).
\end{array}
\]
The universal elements
\[
{\cal R}_{{\cal J}{\cal E}}= \left( {\cal F}_{{\cal J}
{\cal E}}^{{\cal P}}\right) _{21}\left( {\cal F}_{{\cal J} {\cal E}
}^{{\cal P}}\right) ^{-1},\qquad {\cal R}_{{\cal RE}}=
\left( {\cal F}_{{\cal RE}}\right) _{21}
\left( {\cal F}_{{\cal RE}}\right) ^{-1}
\]
are the quantizations of the classical $r$--matrices
\begin{equation}
\label{sl3-pr}
\begin{array}{rcl}
r_{{\cal J} {\cal E}} & = & H^{\cal P}\wedge E+A\wedge B
+\xi H^{\perp}\wedge A, \\[0.2cm]
r_{{\cal RE}} & = & H^{\cal P} \wedge E+A\wedge B+\zeta H^{\perp }
\wedge E .
\end{array}
\end{equation}
They define the dual Lie algebras ${\bf M}_{{\cal J}^{\perp }}^{\ast }$,
${\bf M}_{{\cal R}}^{\ast }$ and the Lie algebra morphisms
${\bf M}^{\ast}\to {\bf M.}$ Applying them to
${\bf M}_{{\cal J}^{\perp}}^{\ast }$ and ${\bf M}_{{\cal R}}^{\ast }$
we get (as an image of this map and after the appropriate
redefinition of the generators)
the 4D subalgebras
${\bf L}_{{\cal J}^{\perp }},{\bf L}_{{\cal R}}\subset $ ${\bf M}$
with  commutators
\[
\begin{array}{rcl}
&  & {\bf L}_{{\cal J}^{\perp }}\left\{
\begin{array}{l}
\lbrack H,E]=E,\quad \lbrack H,A]=A,\quad \lbrack H,B]=0, \\[0.2cm]
\lbrack A,B]=E,\quad \lbrack E,A]=0,\quad \lbrack E,B]=\xi E ;
\end{array}
\right. \\ \\
&  & {\bf L}_{{\cal R}}\ \ \left\{
\begin{array}{l}
\lbrack H,E]=E,\quad \lbrack H,A]=\left( 1-\zeta \right) A,\quad
\lbrack H,B]=\zeta B, \\[0.2cm]
\lbrack A,B]=E,\quad \lbrack E,A]=0,\qquad \qquad\ \ [E,B]=0.
\end{array}
\right.
\end{array}
\]
The algebras ${\bf L}_{{\cal J}^{\perp }},{\bf L}_{{\cal R}}$ and
${\bf L}^{{\cal P}}$ are inequivalent. The deforming functions
$\mu _{{\cal J}^{\perp }}(E,B)=E$ and
$ \left\{\mu _{{\cal R}}(H,B)= B,\; \mu _{{\cal R}}(H,A)= -A \right\}$
are cohomologically nontrivial: $\mu _{{\cal J}^{\perp }},\;
\mu _{{\cal R}}\in H^{2} ( {\bf L}^{{\cal P}},{\bf L}^{{\cal P}} ) $.
Here it is worthy to mention
 that both deforming functions describe
not only the ``tangent vector" to the deformation curve but  the complete
(first order) deformation of ${\bf M}$.
Notice that  ${\bf L}_{{\cal J}^{\perp }}$ and ${\bf L}_{{\cal R}}$
are Frobenius Lie algebras. The
corresponding nondegenerate forms are coboundaries.

In this simple example we have demonstrated that the carrier
subalgebra associated to
the peripheric twist
${\cal F}_{\cal E}^{\cal P}$ enlarged by
${\mit\Phi} _{\cal J}^{\perp }$ or
${\mit\Phi} _{{\cal R}}$
differs nontrivially from the initial ${\bf L}^{{\cal P}}$.

\sect{Construction of peripheric chains}
\label{periphericchains}

The constructions used in the previous Section can be
generalized for the case
where the initial carrier algebra contains the
maximal nilpotent subalgebra
of a simple Lie algebra. The only known twists with such properties
are the (full) chains of extended Jordanian twists (\ref{mainchain}).
In order to use chains of twists for our purposes
we must find their analogs with the
peripheric properties. In \cite{KWL} it was proved that the
peripheric chains exist.

In this section we shall compose and study the full peripheric chains for
${\mathfrak g}=sl(N)$.
We normalize the Cartan elements as $H_{i,k}=(E_{ii}-E_{kk})/2$
and use the
standard $gl(N)$--basis $\{E_{i,j}\}_{i,j=1,\dots N}$.
Consider the canonical chain (\ref{mainchain}) of extended twists
${\cal F}_{{\cal B}_{p\prec 0}}\equiv {\cal F}_{{\cal B}_{p}}
{\cal F}_{{\cal B}_{p-1}}\ldots {\cal F}_{{\cal B}_{0}}$ with
\begin{equation} \label{can-link}
{\cal F}_{{\cal B}_{k}}={\mit\Phi} _{{\cal E}_{k}}
{\mit\Phi} _{{\cal J}_{k}} =
\left( \prod_{s=k+2}^{N-k-1}e^{E_{k+1,s}\otimes
E_{s,N-k}e^{-\frac{1}{2}\sigma _{k+1,N-k}}}\right) e^{H_{k+1,N-k}\otimes
\sigma _{k+1,N-k}},
\end{equation}
where $k=0,\ldots ,p$ and $\sigma_{ij}=\ln (1+E_{i,j})$.
This chain produces the
deformation
${\cal F}_{{\cal B}_{p\prec 0}}:U(sl(N)) \longrightarrow
U_{\cal P} (sl(N))$.
The deformed Hopf algebra
$U_{{\cal B}}(sl(N))$ admits  additional twists.
One of them ${\mit\Phi} _{{\cal J}_{p\prec 0}}^{{\cal R}}$,
a multidimensional analog of
$ ( {\mit\Phi} _{{\cal R}} ) ^{-1}$, must transform $U_{{\cal B}}(sl(N))$
into the quantized algebra of  peripheric type
${\mit\Phi} _{{\cal J}_{p\prec 0}}^{{\cal R}}:U_{{\cal B}}(sl(N))\to
U_{{\cal BP}_{p\prec 0}}(sl(N))$.
To start the construction of the Hopf algebra
$U_{{\cal BP}_{p\prec 0}}(sl(N)) $  we separate the
peripheric factors in the Jordanian twists:
\[
{\mit\Phi} _{{\cal J}_{k}}= e^{H_{k+1,N-k}\otimes \sigma
_{k+1,N-k}}=e^{H_{k+1,N-k}^{{\cal R}}\otimes \sigma
_{k+1,N-k}}e^{H_{k+1,N-k}^{{\cal P}}\otimes \sigma _{k+1,N-k}}=
{\mit\Phi} _{{\cal J}_{k}}^{{\cal R}}{\mit\Phi} _{{\cal J}_{k}}^{{\cal P}}.
\]
Here
\begin{eqnarray}
\label{per-cart}
H_{k+1,N-k}^{{\cal P}} & = &\displaystyle{  \frac{1}{N}
\left(NE_{k+1,k+1}-\sum_{s=1}^{N}E_{s,s}\right)}, \\[0.3cm]
\label{resh-cart}
H_{k+1,N-k}^{{\cal R}} & = & H_{k+1,N-k} - H_{k+1,N-k}^{{\cal P}}.
\end{eqnarray}
The factor ${\mit\Phi} _{{\cal J}_{k}}^{{\cal R}}$ can be
dragged to the end of the extended twist ${\cal F}_{{\cal B}_{k}}$, the
extension factors will be changed:
\begin{equation}
\label{fact-link}
{\cal F}_{{\cal B}_{k}}={\mit\Phi} _{{\cal J}_{k}}^{{\cal R}}
\left(\prod_{s=k+2}^{N-k-1}e^{E_{k+1,s}\otimes E_{s,N-k}}\right)
e^{H_{k+1,N-k}^{{\cal P}}\otimes \sigma _{k+1,N-k}}=
{\mit\Phi} _{{\cal J}_{k}}^{{\cal R}}{\cal  F}_{{\cal B}_{k}}^{{\cal P}}.
\end{equation}
One can check that ${\cal F}_{{\cal B}_{k}}^{{\cal P}}$'s are the
multidimensional peripheric extended twists. In particular, the Cartan
element $H_{k+1,N-k}^{{\cal P}}$ commutes with all the second constituent
generators $ \{ E_{s,N-k}; $ $  s=k+2,\ldots ,N-k-1 \} $ (in the formula
(\ref{gen-ext}) these are the elements
$L_{{\lambda_{0}^{k}} -\lambda ^{\prime }}$).
Applying ${\cal F}_{{\cal B}_{k}}$ to the primitive subalgebra $U(sl(N-2k))$
(this is just the situation that happens after the action
of the $k$ first factors
${\cal F}_{{\cal  B}_{k-1}}\ldots {\cal F}_{{\cal B}_{0}}$) it preserves
the primitivity of all the first constituent generators
$\{ E_{k+1,s}\;;\; s=k+2,\ldots ,N-k-1 \}$.

The factor ${\mit\Phi} _{{\cal J}_{k}}^{{\cal R}}$ can be further pushed to
the very end of the chain because it commutes with all the subsequent links
$\{ {\cal F}_{{\cal B}_{s}};\ s>k \} $. The factors ${\mit\Phi} _{{\cal
J}_{k}}^{{\cal R}}$ with different indices $k$ commute with
each other. Performing the factorizations (\ref{fact-link})
in all the links and collecting the factors
${\mit\Phi} _{{\cal J}_{k}}^{{\cal R}}$ at the end
of the chain we get the following expression
\begin{equation}
\label{chain2}
\begin{array}{rcl}
{\cal F}_{{\cal B}_{p\prec 0}}\equiv {\cal F}_{{\cal B}_{p}}
{\cal F}_{{\cal B}_{p-1}}\ldots {\cal F}_{{\cal B}_{0}} & = &
{\mit\Phi} _{{\cal J}_{p}}^{{\cal R}}{\mit\Phi} _{{\cal J}_{p-1}}^{{\cal R}}
\ldots {\mit\Phi}_{{\cal J}_{0}}^{ {\cal R}}{\cal F}_{{\cal B}_{p}}^{{\cal P}}
{\cal F}_{{\cal B}_{p-1}}^{{\cal P }}\ldots {\cal F}_{{\cal B}_{0}}^{{\cal P}}
\\[0.3cm]
& = & {\mit\Phi} _{{\cal J}_{p\prec 0}}^{{\cal R}}
{\cal F}_{{\cal B}_{p\prec 0}}^{{\cal P}} \ .
\end{array}
\end{equation}

We know \cite{KLO} that ${\cal F}_{{\cal B}_{p\prec 0}}$
is the twist for $U(sl(N))$. In the deformed algebra
$U_{{\cal B}_{p\prec 0}}(sl(N))$ we have
$( p+1 ) $ primitive elements of the type $\sigma _{k+1,N-k}$. It can be
easily checked that all the elements $H_{k+1,N-k}^{{\cal R}}$ are also
primitive in $U_{{\cal B}_{p\prec 0}}(sl(N))$. So, the integral factor
${\mit\Phi}_{ {\cal J}_{p\prec 0}}^{{\cal R}}$ consists of
Reshetikhin twists for commuting elements
$\{ H_{k+1,N-k}^{{\cal R}},\ E_{k+1,N-k};$ $k=0,\ldots ,p\}$.
This means that ${\mit\Phi} _{{\cal J}_{p\prec 0}}^{
{\cal R}}$ is a twist for $U_{{\cal B}_{p\prec 0}}(sl(N))$,
the same is true
for $ ( {\mit\Phi} _{{\cal J}_{p\prec 0}}^{{\cal R}} ) ^{-1}$.
Consequently, the composition
$ ( {\mit\Phi} _{{\cal J}_{p\prec 0}}^{{\cal R}}) ^{-1}
{\cal F}_{{\cal B}_{p\prec 0}}$ is a twist for $U(sl(N))$ and
it immediately follows from (\ref{chain2}) that
${\cal F}_{{\cal B}_{p\prec 0}}^{{\cal P}}$ is also a twist,
${\cal F}_{{\cal B}_{p\prec 0}}^{{\cal P} }:U(sl(N))\to
U_{{\cal BP}_{p\prec 0}}(sl(N))$,
\begin{equation}
{\cal F}_{{\cal B}_{p\prec 0}}^{{\cal P}}
=\prod_{k=0}^{\stackrel{p}{\longleftarrow }}{\cal F}_{{\cal B}
_{k}}^{{\cal P}}=\prod_{k=0}^{\stackrel{p}
{\longleftarrow }}\left( \prod_{s=k+2}^{N-k-1}e^{E_{k+1,s}\otimes
E_{s,N-k}}\right) e^{H_{k+1,N-k}^{{\cal P}}
\otimes \sigma_{k+1,N-k}}.
\end{equation}
This composition of basic twisting factors is
the necessary peripheric analogue
of the chain of extended twists for $sl(N)$,
it is called the peripheric chain
of twists. The chain is full when the number
of links (for ${\bf B^+} (sl(N))$) is maximal: $N/2$
for even and $(N+1)/2$ for odd $N$.

\sect{Peripheric chains as quantization tools}
\label{propertieschains}


Let us consider  $N=2n$ for the even case and $N=2n-1$ for the odd one.
The twisted algebra
$U_{{\cal BP}_{p\prec 0}}(sl(N))$ has $z=(p+1)$ primitive elements
$\sigma_{k+1,N-k}$ and $(N-1-z)$ primitive Cartan generators. Each link
${\cal F}_{{\cal B}_{k}}^{{\cal P}}$ of the chain
${\cal F}_{{\cal B}_{p\prec 0}}^{ {\cal P}}$ is a peripheric extended twist.
After applying it to
$U_{{\cal BP} _{ ( k-1 ) \prec 0}}(sl(N))$ we get $(N-2k-3)$
primitive generators $L_{\lambda ^{k\prime }}$
corresponding to the constituent roots
$\lambda ^{k\prime}$ (such that
$\lambda ^{k\prime} +\lambda^{k\prime \prime }=\lambda _{0}^{k}$).
In our case these are the elements
$ \{ E_{k+1,s}\ |\ s=k+2,\ldots, N-k-1 \}$.
The next link ${\cal F}_{{\cal B}_{k+1}}^{{\cal P}}$ preserves the
primitivity only of one of them, $E_{k+1,N-k-1}$. Thus in
$U_{{\cal BP}_{p\prec 0}}(sl(N))$\ (after having applied $z$ links)
we get $(N-p-2)$
additional primitive generators (with respect to the effect of
the canonical chain).
The following statement will be useful for applications.
\begin{lemma}
The `matreshka' effect \cite{KLO} is valid for peripheric chains.
\end{lemma}
{\em Proof}.
The peripheric link ${\cal F}_{{\cal B}_{k}}^{{\cal P}}$ is the
function of the tensor invariant for the subalgebra
${\mathfrak g}_{\lambda_{0}^{k}}^{\perp }$ and being applied to
$U_{{\cal BP}_{( k-1)\prec 0}}({\mathfrak g}_{\lambda _{0}^{k}}^{\perp })$
such twist  cannot change
its coproducts. On the first step of the inductive process (with
${\cal F}_{ {\cal B}_{0}}^{{\cal P}}$) all the vectors in $\mathfrak g$
 have primitive  coproducts. \hfill $\blacksquare$

In particular, the Hopf algebra $U_{{\cal BP}_{p\prec 0}}(sl(N))$ contains
the subalgebra $U ( sl(N-2z) ) $ whose generators are primitive.

In the following we shall consider the most important case: the  full
peripheric chains ${\cal F}_{{\cal B}_{ ( N-n-1 ) \prec 0}}^{{\cal P}}$
for the enveloping algebra $U(sl(N))$. They have $z=(N-n)$ links
and the twisted algebras $U_{ {\cal BP}_{( N-n-1) \prec 0}}(sl(N))$ contain
$(N-n)$ primitive elements
$\sigma _{k}\equiv \sigma _{k+1,N-k}$, $(n-1)$ primitive Cartan
generators and $(n-1)$ additional primitive generators
$E_{k}^{{\cal P}}\equiv E_{k+1,N-k-1}$. This provides  possibilities
to enlarge  the sequence of twists in the peripheric chain by the new twisting
factors. As it was demonstrated in
Section~\ref{deformationcarriers} these possibilities
are of two kinds: additional Reshetikhin
and additional Jordanian basic twisting factors.

\subsect{Additional Reshetikhin BTF's}

Notice that the primitive elements $I_{m}\in \{ \sigma _{k},
E_{l}^{{\cal P}}|k,l=0,\ldots ,z-1 \} $ commute.
Consequently, the following factor
\[
{\cal F}_{{\cal R}}=\exp \left\{ \beta ^{mn}I_{m}\otimes I_{n}\right\}
\]
is the solution of the Drinfeld equations (\ref{drinf}) for the algebra
 $U_{{\cal BP}_{ ( N-n-1 ) \prec 0}}(sl(N))$. The coefficients
$\beta^{mn}$ are arbitrary. Thus, the composition
\[
{\cal F}_{{\cal RB}}^{{\cal P}}={\cal F}_{{\cal R}}{\cal F}_{{\cal B}
_{\left( N-n-1\right) \prec 0}}^{{\cal P}}
\]
is the twisting element for $U(sl(N))$,
\[
{\cal F}_{{\cal RB}}^{{\cal P}}:U(sl(N))\to U_{{\cal RBP}
_{\left( N-n-1\right) \prec 0}}(sl(N)) .
\]
Each link ${\cal F}_{{\cal B}_{k}}^{{\cal P}}$ in ${\cal F}_{{\cal B}
_{ ( N-n-1 ) \prec 0}}^{{\cal P}}$ can have an independent parameter
$\psi _{k+1}$\ \cite{KLO},
\begin{equation}
 \label{par-link}
{\cal F}_{{\cal B}_{k}}^{{\cal P}}\left( \psi _{k+1}\right) =\left(
\prod_{s=k+2}^{N-k-1}e^{\psi _{k+1}E_{k+1,s}\otimes E_{s,N-k}}\right)
e^{H_{k+1,N-k}^{{\cal P}}\otimes \sigma _{k+1,N-k}\left( \psi _{k+1}\right)
}.
\end{equation}
Let all the variables of ${\cal F}_{{\cal
RB}}^{{\cal P}}$ be proportional to the overall deformation parameter $\xi $,
i.e.
 $\beta ^{mn}\Rightarrow \xi  \beta ^{mn}$ and $\psi _{k+1}\Rightarrow \xi
\psi_{k+1}$. The universal element (\ref{def-t})
\[
{\cal R}_{{\cal RB}}^{{\cal P}}\left( \xi \right)
 =\left( {\cal F}_{{\cal RB}}^{{\cal P}}\left( \xi \right) \right) _{21}
\left( {\cal F}_{{\cal RB}}^{ {\cal P}}\left( \xi \right) \right) ^{-1}
\]
corresponds to the classical $r$--matrix
\begin{equation}
\begin{array}{rcl}
 \label{r-prb}
r_{{\cal RB}}^{{\cal P}} &=&\sum_{k=0}^{p}\psi _{k+1}
\left( H_{k+1,N-k}^{ {\cal P}}\wedge E_{k+1,N-k}
+\sum_{s=k+2}^{N-k-1}E_{k+1,s}\wedge E_{s,N-k}\right)  \\[0.3cm]
&&\qquad +\sum_{n,m=1}^{2z}\beta ^{mn}J_{m}\wedge J_{n}.
\end{array}
\end{equation}
Here $J_{m}\in \{ E_{k+1,N-k},\ E_{l}^{{\cal P}}\ |\ k,l=0,\ldots ,z-1 \} $.
The difference between $r_{{\cal RB}}^{{\cal P}}$ and the
$r$--matrix for the canonical chain
\[
r_{{\cal B}}=\sum_{k=0}^{p}\psi _{k+1}\left( H_{k+1,N-k}\wedge
E_{k+1,N-k} +\sum_{s=k+2}^{N-k-1}E_{k+1,s}\wedge E_{s,N-k}\right)
\]
is essential. This is clearly seen when the corresponding Frobenius forms
are compared
\[
\omega _{{\cal RB}}^{{\cal P}}=\sum_{k=0}^{p}\chi _{k+1}
E_{k+1,N-k}^{\ast}([,])
+\sum_{n,m=1}^{2z}\phi ^{mn}K_{m}^{\ast }\wedge K_{n}^{\ast },
\]
and
\[
\omega _{{\cal B}}=\sum_{k=0}^{p}\chi _{k+1}E_{k+1,N-k}^{\ast }([,]),
\]
where
$K_{m}\in \{ H_{k+1,N-k}^{{\cal P}},\ E_{l}^{{\cal P}}\ |\ k,l=0,\ldots ,z-1 \}
$.
The coefficients $\chi _{k+1}$ and $\phi ^{mn}$ are proportional to
$\psi _{k+1}$
and $\beta ^{mn}$ respectively. Notice that
$\omega _{{\cal RB}}^{{\cal P}}=\omega _{{\cal B} }
+\sum_{n,m=1}^{2z}\phi ^{mn}K_{m}^{\ast}\wedge K_{n}^ {\ast }$
and when all $\phi ^{mn}$ are equal
to zero the forms coincide despite the fact that the
corresponding terms in the $r$--matrix contain different Cartan generators.
Moreover, as we had already mentioned in Section~\ref{periphericchains},
the carrier algebras for
$ r_{{\cal B}}$ and $r_{{\cal RB}}^{{\cal P}} ( \phi ^{mn}=0 ) $
are different. In the general case
(with nonzero $\phi ^{mn}$) the $\omega $--forms are also different: $
\omega _{{\cal B}}$ is a coboundary while\ $\omega _{{\cal RB}}^{{\cal P}}$
is cohomologically nontrivial.

The result is that the peripheric chains provide the possibility to quantize
explicitly the new class of classical $r$--matrices (\ref{r-prb}), defined on
carrier subalgebras of the ${\bf L}^{{\cal P}}$--type, whose  $\omega $--forms
are cohomologically nontrivial.

\subsect{Additional Jordanian BTF's}

Let us consider  the peripheric chain
${\cal F}_{{\cal B}_{( N-n-1) \prec 0}}^{{\cal P}}$ and its
possible extensions
related to the primitive Cartan generators $\{ H_{i}^{\perp }\} $
and the additional primitive elements $\{ E_{i}^{{\cal P}}\} $ ($
i=1,\ldots ,n-1$) belonging to the twisted Hopf algebra $U_{{\cal BP}
_{ ( N-n-1 ) \prec 0}}(sl(N))$. Here, we consider only the
generators $E_{i}^{{\cal P}}$ because the alternative pairs
$ \{ H_{i}^{\perp },\sigma _{k} \} $ appear also in the canonical case and
were treated earlier (see for example \cite{KLS}). The dual elements $
 \{ H_{i}^{\perp \ast } \} $ generate in the root space the
hyperplane orthogonal to the initial roots
$ \{ \lambda_{0}^{k}\;|\;k=0,\ldots ,z-1 \}
$ of the chain ${\cal F}_{{\cal B}_{ ( N-n-1 ) \prec 0}}^{{\cal P}}$. On this
hyperplane we fix the set of basic elements $ \{ H_{i}^{\perp } \} $ and
consider it together with the set $ \{ E_{i}^{{\cal P}} \} $:
\begin{equation}
\label{add-gen}
\begin{array}{rcl}
H_{i}^{\perp } &=&\frac{N-2i}{N}\sum_{l=1}^{N}E_{l,l}-
\frac{2i}{N}\sum_{m=i+1}^{N-i}E_{m,m}, \\[0.25cm]
E_{j}^{{\cal P}} &=&E_{j,N-j},
\end{array}
\end{equation}
with $i,j =1,\ldots ,n-1$. On the joint space generated by
$\{ H_{i}^{\perp},E_{j}^{{\cal P} }\} $ we have a Lie algebra
with the nonzero commutators
\[
\left[ H_{i}^{\perp },E_{j}^{{\cal P}}\right] =\delta _{ij}E_{j}^{{\cal P}}.
\]
In other words, this is the direct sum of $(n-1)$ 2D Borel
subalgebras ${\bf B}_{i}$. This guarantees the possibility to enlarge the
twisting element ${\cal F}_{{\cal B}_{ ( N-n-1 ) \prec 0}}^{{\cal P}}$
with the additional (independent) Jordanian basic factors
${\mit\Phi}_{J_{i}}=e^{H_{i}^{\perp }\otimes \sigma _{i}}$.

Let ${\cal F}_{{\cal JB}}^{{\cal P}}$ be the maximal enlargement of the
full peripheric chain with  $(n-1)$ additional JBF's,
\[
{\cal F}_{{\cal JB}}^{{\cal P}}=\left( \prod_{i=1}^{n-1}e^{H_{i}^{\perp
}\otimes \sigma _{i}}\right) {\cal F}_{{\cal B}_{\left( N-n-1\right) \prec
0}}^{{\cal P}}.
\]
Formally, the carrier subalgebra of this chain coincides with the Borel
subalgebra ${\bf B}^{+} ( sl(N) ) $. It contains all the raising
operators $ \{ E_{lm}\;|\;l,m=1,\ldots ,N\; ;\; l<m \} $, $(N-n)$  generators
$H_{k,N-k}^{{\cal P}}$ and $(n-1)$ generators $H_{i}^{\perp }$. However, as we
shall demonstrate below, the dimension of the carrier for
the classical $r$--matrix
$r_{{\cal JB}}^{{\cal P}}$ remains the same as in the full chain
${\cal F}_{{\cal B}_{ ( N-n-1 ) \prec 0}}^{{\cal P}}$.

To study the properties of the enlarged chain in detail we must introduce
in ${\cal F}_{{\cal JB}}^{{\cal P}}$ the full set of parameters. This can be
done by means of an involutive automorphism of the carrier subalgebra and
the discrete transformations of subalgebras inside the carrier. It is
sufficient to consider the maximal nilpotent subalgebra
${\bf N}^{+} (sl(N) ) \subset {\bf B}^{+} ( sl(N) ) $. We perform the smooth
scaling of the generators in ${\bf N}^{+} ( sl(N) ) $:
\[
\left\{ E_{lm}\Rightarrow \alpha _{lm}E_{lm}\;|\; \alpha _{lm}\in {\C}
;\; l,m=1,\ldots ,N;\; l<m\right\} .
\]
The scaling factors $\alpha _{lm}$ are as follows:
\begin{equation}\label{alphas}
\begin{array}{lll}
\textrm{for} &\quad \left\{ E_{lm}\left|
\begin{array}{c}
l=1,\ldots ,z-1, \\[0.2cm]
m=2,\ldots ,z, \\[0.2cm]
l<m
\end{array}
\right. \right\} & \alpha _{lm}=\displaystyle \frac{\psi _{l}}{\psi _{m}}, \\
\\
\textrm{for} &\quad \left\{ E_{lm}\left|
\begin{array}{c}
l=1,\ldots ,z, \\[0.2cm]
m=z+1,\ldots ,N,
\end{array}
\right. \right\} & \alpha _{lm}=\psi _{l}\zeta _{\left( N-m\right) }, \\
\\
\textrm{for} &\quad \left\{ E_{lm}\left|
\begin{array}{c}
l=z+1,\ldots ,N-1, \\[0.2cm]
m=z+2,\ldots ,N, \\[0.2cm]
l<m
\end{array}
\right. \right\} & \alpha _{lm}=\displaystyle \frac{\zeta _{\left( N-m\right) }
}
{\zeta _{\left( N-l\right) }},
\end{array}
\zeta _{0} = 1.
\end{equation}
All the parameters $ \{ \psi _{l},\zeta _{i}\;|\;l=1,\ldots,z;\;
i=1,\ldots ,N-z-1 \} $ are independent. The first subset $ \{\psi _{l}\} $
consists of the parameters of  links, they were
already introduced in ${\cal F}_{{\cal B}_{k}}^{{\cal P}} ( \psi_{k+1} ) $
(see (\ref{par-link})), the entries of the second subset
$\{ \zeta _{i} \} $ refer to the Jordanian factors
${\mit\Phi}_{J_{i}}=e^{H_{i}^{\perp }\otimes \sigma _{i}}$.
Notice that each argument
$E_{i}^{{\cal P}}$ of $\sigma _{i}$\ is  scaled by a product of
parameters, $E_{i}^{{\cal P}}\Rightarrow  ( \psi _{i}\zeta
_{j} ) E_{i}^{{\cal P}}$, because it is already scaled in the
corresponding link ${\cal F}_{{\cal B}_{i-1}}^{{\cal P}} ( \psi _{i} ) $
of the chain.
The discrete parameters describe the property: in any link of the chain the
extension factor can be switched off. Only the full extension
\[
{\mit\Phi} _{{\cal E}_{k}}\left( \psi _{k+1}\kappa _{k+1} \right)
=\prod_{\lambda ^{\prime }\in \pi _{k}^{\prime }}{\mit\Phi} _{{\cal E}
_{\lambda ^{\prime }}}\left( \psi _{k+1}\kappa _{k+1}\right)
=\prod_{s=k+2}^{N-k-1}e^{\psi _{k+1}\kappa _{k+1}E_{k+1,s}\otimes E_{s,N-k}}
\]
has the character of the basic extending factor (\ref{ext}) with the continuous
parameter $\psi _{k+1}$ and the discrete parameter $\kappa _{k+1}=0,1$.
The separate EF's inside ${\mit\Phi}_{{\cal E}_{k}}$ conserve their
independence and can be switched off. But, in the general case, such
cancellation ruins the matreshka effect and the structure of the chain is,
 hence, lost. The parameterized enlarged chain has the form
\begin{equation}\label{par-en-chain}
\begin{array}{l}\hskip -0.4cm
{\cal F}_{{\cal JB}}^{{\cal P}}\left( \left\{ \psi _{l},
\kappa_{l},\zeta _{i}\right\} \right) \\[0.25cm]
\hskip -0.3cm
=\left( \prod\limits_{i=1}^{n-1}e^{H_{i}^{\perp }
\otimes \sigma_{i}\left(\psi _{i}\kappa _{i}\zeta _{i}\right) }\right)
{\cal F}_{{\cal B} _{\left( N-n-1\right) \prec 0}}^{{\cal P}}
\left( \left\{ \psi _{l},\kappa_{l},\zeta _{i}\right\} \right)
\\[0.25cm] \hskip -0.3cm
=\left( \prod\limits_{i=1}^{n-1}e^{H_{i}^{\perp }\otimes \sigma_{i}
\left( \psi _{i}\kappa _{i}\zeta _{i}\right) }\right)
\prod\limits_{k=0}^{\stackrel{z}{\longleftarrow }}
\left( \prod\limits_{s=k+2}^{N-k-1}e^{\psi_{k+1}\zeta _{k}
\kappa _{k+1}E_{k+1,s}\otimes E_{s,N-k}}\right) e^{H_{k+1,N-k}^{{\cal P}}
\otimes \sigma_{k+1,N-k}\left( \psi _{k+1}\zeta _{k}\right) }.
\end{array}
\end{equation}
Any number of links ${\cal F}_{{\cal B}_{i-1}}^{{\cal P}}
(  \{ \psi_{l},\kappa _{l},\zeta _{i} \}  ) $ and JBF's ${\mit\Phi}
_{J_{i}} ( \psi _{i}\kappa _{i}\zeta _{i} ) =e^{H_{i}^{\perp}
\otimes \sigma _{i} ( \psi _{i}\kappa _{i}\zeta _{i} ) }$ can be
switched off in the enlarged chain. This is ensured by the following
rearrangement of parameters:
\[
\psi _{l}=\nu _{l}\prod_{r=1}^{l-1}\frac{\nu _{r}}{\rho _{r}},\qquad
\zeta _{i}=\prod_{r=1}^{i}\frac{\rho _{r}}{\nu _{r}}.
\]
\begin{equation}\label{par-en-2}
\begin{array}{l}\hskip -0.3cm
{\cal F}_{{\cal JB}}^{{\cal P}}\left( \left\{ \nu _{l},\kappa _{l},
\rho_{i}\right\} \right)  \\[0.25cm] \hskip -0.2cm
\quad =\left( \prod\limits_{i=1}^{n-1}e^{H_{i}^{\perp }\otimes
\sigma _{i}\left(\kappa _{i}\rho _{i}\right) }\right)
{\cal F}_{{\cal B}_{\left( N-n-1\right) \prec 0}}^{{\cal P}}
\left( \left\{ \nu _{l},\kappa_{l}\right\} \right)
\\[0.25cm] \hskip -0.2cm
\quad =\left( \prod\limits_{i=1}^{n-1}e^{H_{i}^{\perp }
\otimes \sigma _{i}\left( \kappa _{i}\rho _{i}\right) }\right)
\prod\limits_{k=0}^{\stackrel{z}
{\longleftarrow }}\left( \prod\limits_{s=k+2}^{N-k-1}e^{\nu _{k+1}
\kappa_{k+1}E_{k+1,s}\otimes E_{s,N-k}}\right)
e^{H_{k+1,N-k}^{{\cal P}}\otimes
\sigma _{k+1,N-k}\left( \nu _{k+1}\right) }.
\end{array}
\end{equation}
From the structural point of view the independence of the
additional factors
${\mit\Phi} _{J_{i}} ( \kappa _{i}\rho _{i} ) $ is obvious.
Due to the matreshka effect switching off a link
${\cal F}_{{\cal B}_{ ( N-n-1 ) \prec 0}}^{{\cal P}}
(  \{ \nu _{l},\kappa _{l} \}  ) $ in
${\cal F}_{{\cal JB}}^{{\cal P}} (  \{ \nu _{l},\rho _{i},
\kappa_{l} \}  ) $   can not prevent any of $E_{i}^{{\cal P}}$'s
to be primitive in
$U_{{\cal B} ( N-n-1 ) \prec 0}^{{\cal P}}(sl(N))$. When an
extension ${\mit\Phi} _{{\cal E}_{l}} ( \nu _{l}\kappa _{l} ) $ is
switched off by putting $\kappa _{l}=0$ the corresponding factor
$e^{H_{l}^{\perp }\otimes \sigma _{l} ( \kappa _{l}\rho _{l} ) }$
also vanishes. This correlates with the fact that when the extension
${\mit\Phi}_{{\cal E}_{l}} ( \nu _{l}\kappa _{l} ) $ is absent (the
corresponding link contains only a JF) the coproduct
$\Delta _{{\cal J} } ( E_{l}^{{\cal P}} ) $
is not primitive  and the additional JF
${\mit\Phi} _{J_{l}}$ has to be canceled.
Strictly speaking, even  in such a
situation the number of additional factors ${\mit\Phi} _{J_{i}} $
can be preserved.
Instead of $\Delta _{{\cal J}} ( E_{l}^{{\cal P}} ) $
we have the primitive
$\Delta _{{\cal J}} ( E_{N-l,N-l+1} ) $ and the factor $
e^{H_{l}^{\perp }\otimes \sigma _{l} (\kappa _{l}\rho _{l} ) }$
can be substituted by
$e^{H_{l}^{\perp }\otimes \sigma _{N-l,N-l+1} (\rho _{l} ) }$.

When the extensions are switched off in all the links
the initial chain degenerates into the twist
\[
{\cal F}_{{\cal B}}^{{\cal P}}\left( \left\{ \nu _{l},0,\rho _{i}\right\}
\right) =\prod\limits_{k=0}^{z}e^{H_{k+1,N-k}^{{\cal P}}\otimes
\sigma_{k+1,N-k}\left( \nu _{k+1}\right) }
\]
(with the special
``peripheric'' choice of the Cartan generators $H_{k+1,N-k}^{{\cal P}}$).
This multi-Jordanian twist can be enlarged by $(n-1)$ JF's,
\[
{\cal F}_{{\cal JB}}^{{\cal P}}\left( \left\{ \nu _{l},0,\rho _{i}\right\}
\right) =\left( \prod\limits_{i=1}^{n-1}e^{H_{i}^{\perp }\otimes
\sigma_{N-l,N-l+1}\left( \rho _{i}\right) }\right)
\prod\limits_{k=0}^{z}e^{H_{k+1,N-k}^{{\cal P}}\otimes \sigma
_{k+1,N-k}\left( \nu _{k+1}\right) }.
\]

Now we shall consider the quasi-classical limit for the enlarged chain
${\cal F}_{{\cal JB}}^{{\cal P}} ( \{ \nu _{l},\kappa _{l},\rho_{i} \}) $.
To simplify the formulas we put
$ \{ \kappa _{l}=1\;|\;l=1,\ldots ,n-1 \} $. The quasi-classical
limit for ${\cal F}_{{\cal JB}}^{{\cal P}} ( \{ \nu _{l},1,\rho
_{i} \}  ) \equiv {\cal F}_{{\cal JB}}^{{\cal P}} (  \{
\nu _{l},\rho _{i} \}  ) $ is obtained through the
substitution $\psi _{k}\Rightarrow \xi \psi _{k}$ in (\ref{par-en-chain}) or
$\nu _{l}\Rightarrow \xi \nu _{l},\ \rho _{i}\Rightarrow \xi \rho _{i}$ in
(\ref{par-en-2}) with the overall deformation parameter $\xi $. In the
neighborhood of the origin the ${\cal R}$--matrix
${\cal R}_{{\cal JB}}^{{\cal P}}$ has the expansion
\begin{eqnarray*}
{\cal R}_{{\cal JB}}^{{\cal P}}\left( \left\{ \psi _{l},\zeta _{i};\xi
\right\} \right) &=&1\otimes 1+\xi r_{{\cal JB}}^{{\cal P}}\left( \left\{
\psi _{l},\zeta _{i}\right\} \right) +{\cal O}\left( \xi \right)  \\[0.2cm]
&=&1\otimes 1+\xi r_{{\cal JB}}^{{\cal P}}\left( \left\{ \nu _{l},\rho
_{i}\right\} \right) +{\cal O}\left( \xi \right) .
\end{eqnarray*}
This means that the deformation $U(sl(N))\to U_{{\cal JB}\left(
N-n-1\right) \prec 0}^{{\cal P}}(sl(N))$ performed by the twist
${\cal F}_{{\cal JB}}^{{\cal P}}\left( \left\{ \psi _{l},\zeta _{i};\xi
\right\} \right) $
can be treated as a quantization of the classical mechanical system
described by the $r$--matrix
\begin{equation}\label{r-jordchain}
\begin{array}{lll}
r_{{\cal JB}}^{{\cal P}} & = & \displaystyle
{\sum_{k=0}^{p}\psi _{k+1}\zeta _{k}
\left( H_{k+1,N-k}^{{\cal P}}\wedge E_{k+1,N-k}+\sum_{s=k+2}^{N-k-1}
E_{k+1,s}\wedge E_{s,N-k}\right) }\\[0.2cm]
&  &\qquad  +\displaystyle{\sum_{i=1}^{n-1}\psi _{i}\zeta _{i}H_{i}^{\perp }
\wedge E_{i,N-i} } \\[0.3cm]
& =& \displaystyle{\sum_{k=0}^{p}\nu _{k+1}\left( H_{k+1,N-k}^{{\cal P}}
\wedge E_{k+1,N-k}+\sum_{s=k+2}^{N-k-1}E_{k+1,s}\wedge E_{s,N-k}\right) }
\\[0.2cm]
&  & \qquad \displaystyle{+\sum_{i=1}^{n-1}\rho _{i}H_{i}^{\perp }
\wedge E_{i,N-i}}.
\end{array}
\end{equation}
It can be easily checked that on the space of ${\bf B}^{+} (sl(N) ) $
this $r$--matrix is degenerate.

Now we shall show that the carrier subalgebra of $r_{{\cal JB}}^{{\cal P}}$
has the dimension $  ( N^{2}+N-2n )/2 $ (the same as in the
case of $r_{{\cal B}}^{{\cal P}}$).

\begin{lemma}
The carrier algebra of the classical $r$--matrix
$r_{{\cal JB}}^{{\cal P}} (  \{ \psi _{l},\zeta _{i} \}  )$
(\ref{r-jordchain}) is equivalent to the
$ (( N^{2}+N-2n)/2)$--dimensional subalgebra of
${\bf B}^{+} (sl(N) )$ generated by the following sets of elements:
\begin{equation}
\label{car-alg}
\begin{array}{ll}
H_{i,N-i+1}^{{\cal P}} &\qquad \ i=1,\ldots ,z \\[0.25cm]
H_{j}^{\perp } & \qquad \ j=1,\ldots ,n-1 \\[0.25cm]
E_{lm} &\quad
\left\{
\begin{array}{l}
l=1,\ldots ,z-1\ ;\ \ m=2,\ldots ,z\ ;\ \ l<m \\[0.2cm]
l=1,\ldots ,z\ ;\ \ m=z+1,\ldots ,N \\[0.2cm]
l=z+1,\ldots ,N-2\ ;\ \ m=z+3,\ldots ,N\ ;\ \ l<m
\end{array} \right\}
\end{array}
\end{equation}
\end{lemma}
{\em Proof}. It is obvious that the set (\ref{car-alg})
generates a subalgebra, that we denote
${\mathfrak g}_{{\cal JB}}^{{\cal P}}$. It
contains the Cartan subalgebra of $sl(N)$ and all the positive roots
operators except those corresponding to the following subset of basic roots
\[
\left\{ \lambda _{z+i}=e_{z+i}-e_{z+i+1}\;|\; i=1,\ldots ,N-z-1\right\} .
\]

Let us construct the smooth set of injections $\varphi $ of
${\mathfrak g}_{{\cal JB}}^{{\cal P}}$ into ${\bf B}^{+} ( sl(N) )$
depending on the
parameters $ \{ \zeta _{i} \} $:
\begin{equation} \label{iject}\hskip-0.3cm
\varphi \left( \left\{ \zeta _{i}\right\} \right) :\left\{
\begin{array}{ll}
H_{i,N-i+1}^{{\cal P}}\rightarrow H_{i,N-i+1}^{{\cal P}} & \ \ i
=1,\ldots ,z\ ; \\[0.2cm]
H_{j}^{\perp }\rightarrow -B_{j}\left( \left\{ \zeta _{i}\right\} \right)
&\ \ j=1,\ldots ,n-1\ ; \\[0.2cm]
E_{lm}\rightarrow E_{lm} &\ \
l=1,\ldots ,z-1\ ;\ \ m=2,\ldots ,z\ ; \ \ l<m\ ; \\[0.2cm]
E_{lm}\rightarrow A_{lm}\left( \left\{ \zeta _{i}\right\} \right) &\ \
l=1,\ldots ,z\ ;\ \ m=z+1,\ldots ,N\ ; \\[0.2cm]
E_{lm}\rightarrow C_{lm}\left( \left\{ \zeta _{i}\right\} \right) &\ \
l=z+1,\ldots ,N-2\ ;\ \ m=z+3,\ldots ,N\ ;\ \ l<m
\end{array}
\right.
\end{equation}
where
\begin{equation}\label{b-gen}
\begin{array}{lll}
B_{j}\left( \left\{ \zeta _{i}\right\} \right) &=&\displaystyle{
\sum_{s=1}^{j}\frac{ \zeta _{j-s}}{\zeta _{j}}E_{N-j,N-j+s}-H_{j}^{\perp }},
\\[0.25cm]
A_{lm}\left( \left\{ \zeta _{i}\right\} \right) &=&\displaystyle{
\frac{1}{\zeta_{N-m}}\sum_{s=0}^{N-m}\zeta _{N-m-s}E_{l,m+s}}, \\[0.25cm]
C_{lm}\left( \left\{ \zeta _{i}\right\} \right)&=&\displaystyle{
E_{lm}-\frac{\zeta _{N-l}}{\zeta _{N-m}}\left(
\sum_{s=m}^{N}\frac{\zeta _{N-s}}{\zeta _{N-l+1}}E_{l-1,s}
+\sum_{s=1}^{N-m}\frac{\zeta _{N-m-s}}{\zeta _{N-l}}E_{l,m+s}\right) }.
\end{array}\end{equation}
Direct computations show that the map
$\varphi  (  \{ \zeta_{i} \}  ) $ is an automorphism of
${\mathfrak g}_{{\cal JB}}^{{\cal P}}$. Let
us denote by $\left\{ D_{s}\right\} $ the basis of
${\mathfrak g}_{{\cal JB}}^{{\cal P}}$
obtained as the image of the natural basis (\ref{car-alg}),
$$
\varphi \left( \left\{ \varsigma _{i}\right\} \right) :\left\{
H_{i,N-i+1}^{{\cal P}},H_{j}^{\perp },E_{pt}\right\} \longrightarrow
\left\{D_{s}\right\} \equiv \left\{ H_{i,N-i+1}^{{\cal P}},
B_{j}\left( \left\{
\varsigma _{i}\right\} \right) ,E_{lm},A_{lm}
\left( \left\{ \varsigma _{i}\right\} \right) ,
C_{lm}\left( \left\{ \varsigma _{i}\right\}
\right)\right\}.
$$
The $r$--matrix $r_{{\cal JB}}^{{\cal P}}\left( \left\{ \psi _{l},
\varsigma_{i}\right\} \right) $
can be expressed in terms of the generators
$\left\{D_{s}\right\} $. First let us rewrite it as follows:
\begin{equation}
\label{r-jpb-new}
\begin{array}{lll}
r_{{\cal JB}}^{{\cal P}}\left( \left\{ \psi _{l},\varsigma _{i}\right\}
\right)
& = & \sum\limits_{k=0}^{p}\psi _{k+1}\varsigma _{k}E_{k+1,N-k-1}\wedge
\left( E_{N-k-1,N-k}-\frac{\varsigma _{k+1}}{\varsigma _{k}}
H_{j+1}^{\perp}\right) \\
&  & +\sum\limits_{k=0}^{p}\psi _{k+1}\varsigma _{k}
\left( H_{k+1,N-k}^{{\cal P}}\wedge E_{k+1,N-k}
+\sum\limits_{s=k+2}^{N-k-2}E_{k+1,s}\wedge E_{s,N-k}\right) .
\end{array}
\end{equation}

Consider the restriction of the map $\varphi _{E}$ (\ref{iject}) to the
subspace generated by
\[
\begin{array}{lll}
E_{lm} &\quad \textrm{for}
&\quad l=1,\ldots ,z-1\ ;\ \ m=2,\ldots ,z\ ;\ \ l<m\ ; \\[0.2cm]
E_{lm} &\quad \textrm{for} & \quad  l=1,\ldots ,z\ ;
\ \ m=z+1,\ldots ,N\ ;\\[0.2cm]
E_{lm} &\quad \textrm{for} &\quad l=z+1,\ldots ,N-2\ ;\ \
m=z+3,\ldots ,N\ ; \ \ l<m\ .
\end{array}
\]
The formulas (\ref{b-gen}) describe the
decompositions $\varphi _{E}\left( E_{ij}\right) =D_{s\left( ij\right)
}=\left( \varphi _{E}\right) _{s}^{pt}E_{pt}$. It is easy to see that the
matrix $\left\{ \left( \varphi _{E}\right) _{s}^{pt}\right\} $ is invertible
and $E_{pt}=$ $\left( \varphi _{E}^{-1}\right) _{pt}^{s}D_{s}$ are the
linear combinations of the elements $E_{ij}$ and $A_{ij}$\ whose indices are
$\left\{ i\geq l,\; j\leq m\right\} $. The remaining entries of
$r_{{\cal JB}}^{{\cal P}}$ (\ref{r-jpb-new}) are the combinations
$(C_{N-k-1,N-k}-\frac{\varsigma _{k+1}}{\varsigma _{k}}H_{j+1}^{\perp })$.
According to the
definition (\ref{b-gen}) these terms depend only on $B_{j}$'s and the
generators $\left\{ C_{pt}\; |\; p=z+1,\ldots ,N-2;\;
t=z+3,\ldots ,N;\; p<t\right\}
$. As a result all the arguments of the $r$--matrix
$r_{{\cal JB}}^{{\cal P}}\left( \left\{ \psi _{l},
\varsigma _{i}\right\} \right) $
 are proved to belong to the space generated by
$H_{i,N-i+1}^{{\cal P}},B_{j},E_{lm},A_{lm}$ and $C_{lm}$,
that is the algebra
${\mathfrak g}_{{\cal JB}}^{{\cal P}}$ described above:
$$
r_{{\cal JB}}^{{\cal P}}\left( \left\{ \psi
_{l},\varsigma _{i}\right\} \right) \in \varphi
\left( {\mathfrak g}_{{\cal JB}}^{{\cal P}}\right) \wedge \varphi
\left( {\mathfrak g}_{{\cal JB}}^{{\cal P}}\right)
\approx {\mathfrak g}_{{\cal JB}}^{{\cal P}}\wedge
{\mathfrak g}_{{\cal JB}}^{{\cal P}} .
$$
In this presentation the $r_{{\cal JB}}^{{\cal P}}$--matrix
is nondegenerate.
This proves the Lemma.\hfill $\blacksquare$

We can calculate the symplectic form
$\omega _{{\cal JB}}^{{\cal P}}$ corresponding to
$r_{\cal JB}^{\cal P} (  \{ \psi _{l},\zeta _{i} \} ) $.
It  is a special kind
of  coboundary generated by the basic forms
$A_{l,m}^{\ast } (  [ \ ,\ ]  )$
with $ \{ l=1,\ldots ,z;\; m=z+1,\ldots ,N-l+1 \}$.
As far as the map $\varphi $ is
an isomorphism we can use the canonical basis (\ref{car-alg}). In this case
the corresponding coboundaries are
$E_{l,m}^{\ast }\left( \left[ ,\right] \right) $. The symplectic form
$\omega _{{\cal JB}}^{{\cal P}}$ has a very simple structure:
\begin{equation}
\label{omega}
\omega _{{\cal JB}}^{{\cal P}}\left(
\left\{ \varsigma _{i}\right\} \right)
=-\sum_{l=1}^{z}\sum_{m=z+1}^{N-l+1}
\frac{1}{\alpha _{lm}}E_{l,m}^{\ast}
\left( \left[ ,\right] \right) ,
\end{equation}
here
$\alpha _{lm}=\alpha _{lm}\left(
\left\{ \varsigma _{i}\right\} \right)$
are the scaling factors defined in (\ref{alphas}).

\sect{Examples}
\label{examples}

We shall demonstrate the properties of the enlarged peripheric chains
presenting explicitly the expressions corresponding to
three special cases:
$U(sl(4))$, $U(sl(7))$ and $U(sl(3))$.

\subsect{$U(sl(4))$}

This is the simplest case where the peripheric properties can
be visualized because four is the\ lowest dimension of the space whose
algebra of linear transformations has the nontrivial chain
of extended twists
\cite{KLO} characterized by the indices
$j=1,\;i=1,2;\; n=2;\; z=N-n=2$,
\begin{equation}
\begin{array}{lll}
{\cal F}_{{\cal B}_{1\prec 0}}^{{\cal P}} &=
&\prod\limits_{k=0}^{\stackrel{1}
{\longleftarrow }}{\cal F}_{{\cal B}_{k}}^{{\cal P}}=
\prod\limits_{k=0}^{\stackrel{1}{\longleftarrow }}\left(
\prod\limits_{s=k+2}^{3-k}e^{E_{k+1,s}\otimes
E_{s,4-k}}\right) e^{H_{k+1,4-k}^{{\cal P}}\otimes \sigma _{k+1,4-k}}
\label{4-chain} \\[6mm]
&=&e^{H_{2,3}^{{\cal P}}\otimes \sigma _{2,3}}e^{E_{1,2}\otimes
E_{2,4}}e^{E_{1,3}\otimes E_{3,4}}
e^{H_{1,4}^{{\cal P}}\otimes \sigma_{1,4}}.
\end{array}\end{equation}
Notice that the extension factor in the second link is trivial.
This chain differs from the canonical one by the Cartan generators in
the Jordanian factors:
\[
\begin{array}{lll}
H_{1,4}^{{\cal P}} & = & \frac{1}{4}\left(
3E_{1,1}-\sum\limits_{s=2}^{4}E_{s,s}\right) , \\[0.35cm]
H_{2,3}^{{\cal P}} & = & \frac{1}{4}\left(
-E_{1,1}+3E_{2,2}-E_{3,3}-E_{4,4}\right) .
\end{array}
\]
The basic vector $ ( H_{1}^{\perp } ) ^{\ast }$ of the hyperplane
orthogonal to
$  \lambda _{1,4}=e_{1}-e_{4}$ and $\lambda_{2,3}=e_{2}-e_{3} $
corresponds to the generator
\[
H_{1}^{\perp }=\frac{1}{2}\left( E_{1,1}-E_{2,2}-E_{3,3}+E_{4,4}\right) .
\]
The coproducts of the generators of ${\bf B}^{+} ( sl(4) ) $ can be
considered as a modification of the coproducts
$\Delta _{_{{\cal B}_{1\prec 0}}}$
(twisted by a canonical chain) due to the Reshetikhin ``rotation'' of
the type $ ( {\mit\Phi} _{{\cal J}_{1\prec 0}}^{{\cal R}} ) ^{-1}$
(see the expression (\ref{chain2})):
\begin{equation}
\label{delta-pb}
\begin{array}{lll}
\Delta _{{\cal B}_{1\prec 0}}^{{\cal P}}\left( E_{1,2}\right)  & = &
E_{1,2}\otimes e^{-\sigma _{1,4}}+1\otimes E_{1,2}-H_{2,3}^{{\cal P}}
\otimes E_{13}e^{-\sigma _{2,3}}; \\[0.25cm]
\Delta _{{\cal B}_{1\prec 0}}^{{\cal P}}\left( E_{1,3}\right)  & = &
E_{1,3}\otimes 1+1\otimes E_{1,3}; \\[0.25cm]
\Delta _{{\cal B}_{1\prec 0}}^{{\cal P}}\left( E_{1,4}\right)  & = &
E_{1,4}\otimes e^{\sigma _{1,4}}+1\otimes E_{1,4}; \\[0.25cm]
\Delta _{{\cal B}_{1\prec 0}}^{{\cal P}}\left( E_{2,3}\right)  & = &
E_{2,3}\otimes e^{\sigma _{2,3}}+1\otimes E_{2,3}; \\[0.25cm]
\Delta _{{\cal B}_{1\prec 0}}^{{\cal P}}\left( E_{2,4}\right)  & = &
E_{2,4}\otimes e^{\sigma _{2,3}}+e^{\sigma _{1,4}}\otimes E_{2,4};
\\[0.25cm]
\Delta _{{\cal B}_{1\prec 0}}^{{\cal P}}\left( E_{3,4}\right)  & = &
E_{3,4}\otimes 1+e^{\sigma _{1,4}}\otimes E_{3,4}
+H_{2,3}^{{\cal P}}\otimes
E_{24}e^{-\sigma _{2,3}}; \\[0.25cm]
\Delta _{{\cal B}_{1\prec 0}}^{{\cal P}}
\left( H_{1,4}^{{\cal P}}\right)
& = & H_{1,4}^{{\cal P}}\otimes e^{-\sigma _{1,4}}
+1\otimes H_{1,4}^{{\cal P}
}-E_{1,3}\otimes E_{3,4}e^{-\sigma _{1,4}} \\[0.2cm]
&& \qquad -\left( E_{1,2}+H_{2,3}^{{\cal P}}E_{1,3}\right)
\otimes E_{2,4}e^{-\sigma_{1,4}-\sigma _{2,3}};\\[0.25cm]
\Delta _{{\cal B}_{1\prec 0}}^{{\cal P}}
\left( H_{2,3}^{{\cal P}}\right)
& = & H_{2,3}^{{\cal P}}\otimes e^{-\sigma _{2,3}}
+1\otimes H_{2,3}^{{\cal P}
}; \\[0.25cm]
\Delta _{{\cal B}_{1\prec 0}}^{{\cal P}}
\left( H_{1}^{\perp }\right)
& = & H_{1}^{\perp }\otimes 1+1\otimes H_{1}^{\perp } .
\end{array}
\end{equation}
We have here two additional primitive elements
\[
E_{1}^{{\cal P}}=E_{1,3},\qquad H_{1}^{\perp },
\]
that determine the Borel algebra
\[
\left[ H_{1}^{\perp },E_{1}^{{\cal P}}\right] =E_{1}^{{\cal P}}.
\]
The twist ${\cal F}_{{\cal B}_{1\prec 0}}^{{\cal P}}$ (\ref{4-chain})
can be enlarged by the factor ${\mit\Phi} _{J_{1}}
=e^{H_{1}^{\perp }\otimes
\sigma_{1}}$, with $\sigma _{1}=\ln  ( 1+E_{1}^{{\cal P}} ) $
\begin{equation}
\label{tw-pjb}
\begin {array}{rcl}
{\cal F}_{{\cal JB}}^{{\cal P}} &=&e^{H_{1}^{\perp }\otimes \sigma _{1}}
{\cal F}_{{\cal B}_{1\prec 0}}^{{\cal P}} \\[0.2cm]
&=&e^{H_{1}^{\perp }\otimes \sigma _{1}}e^{H_{2,3}^{{\cal P}}\otimes
\sigma_{2,3}}e^{E_{1,2}\otimes E_{2,4}}e^{E_{1,3}
\otimes E_{3,4}}e^{H_{1,4}^{{\cal P}}\otimes \sigma _{1,4}}.
\end{array}
\end{equation}
The enlarged twist corresponds to the ${\cal R}$--matrix
\[
{\cal R}_{{\cal JB}}^{{\cal P}}=e^{\sigma _{1}\otimes H_{1}^{\perp }}
\left({\cal F}_{{\cal B}_{1\prec 0}}^{{\cal P}}\right)_{21}
\left( {\cal F}_{{\cal B}_{1\prec 0}}^{{\cal P}}\right)^{-1}
e^{-H_{1}^{\perp }\otimes \sigma _{1}}.
\]
This universal element (supplied with the
full set of deformation parameters
$ \{ \psi _{l}\Rightarrow \xi \psi _{l},\;
\varsigma _{1}\;|\;l=1,2 \}$
with the overall parameter $\xi $) can be considered as
a quantization of the classical $r$--matrix
$r_{{\cal JB}}^{{\cal P}} (  \{ \psi_{l},\varsigma_{i} \}  )$
\begin{equation}
\label{trans-r-mat-4}
\begin{array}{lll}
r_{{\cal JB}}^{{\cal P}}\left( \left\{ \psi _{l},\varsigma _{1}\right\}
\right) & = & \psi _{1}E_{1,3}\wedge \left( E_{3,4}-\varsigma
_{1}H_{1}^{\perp }\right) +\psi _{2}\varsigma _{1}H_{2,3}^{{\cal P}}
\wedge E_{2,3} \\[0.25cm]
&  &\qquad +\; \psi _{1}\left( H_{1,4}^{{\cal P}}\wedge E_{1,4}+E_{1,2}
\wedge E_{2,4}\right) .
\end{array}
\end{equation}
The Lie-Poisson structure fixed by
$r_{{\cal JB}}^{{\cal P}}\left( \left\{ \psi _{l},
\varsigma_{i}\right\} \right) $
can be redefined in terms of the algebra
${\mathfrak g}_{{\cal JB}}^{{\cal P}}$ that in our case is 8D generated by
$\left\{H_{1,4}^{{\cal P}},H_{2,3}^{{\cal P}},
H_{1}^{\perp},E_{p,t}\;|\;p=1,2;\;t=2,3,4;\;p<t\right\} $.
This is an algebra of motion over
the 4D space with the translations
${\mathfrak g}_{{\cal P}}=\left\{E_{p,t}\;|\;p=1,2;\;t=3,4\right\} $
and the subalgebra
${\mathfrak g}_{{\cal H}}$
containing the operator $E_{1,2}$ and the Cartan generators
$\left\{ H_{1,4}^{{\cal P}},\; H_{2,3}^{{\cal P}},\; H_{1}^{\perp }\right\} $,
\[
{\mathfrak g}_{{\cal JB}}^{{\cal P}}=
{\mathfrak g}_{{\cal H}}\vdash {\mathfrak g}_{{\cal P}}.
\]
In terms of the image $\varphi
\left({\mathfrak g}_{{\cal JB}}^{{\cal P}}\right) $
generated by
\[
\begin{array}{lllll}
B_{1}\left( \left\{ \varsigma _{1}\right\} \right) & = &
\frac{1}{\varsigma_{1}}E_{3,4}-H_{1}^{\perp },  & \quad   E_{12}, &\\[0.3cm]
A_{13}\left( \varsigma _{1}\right) & = & E_{1,3}+\frac{1}{\varsigma _{1}}
E_{1,4}, &  \quad  A_{14}=E_{14}, &\quad  H_{1,4}^{{\cal P}}, \\[0.3cm]
A_{23}\left( \varsigma _{1}\right) & = & E_{2,3}+\frac{1}{\varsigma _{1}}
E_{2,4}, & \quad A_{24}=E_{24}, &\quad  H_{2,3}^{{\cal P}},
\end{array}
\]
this $r$--matrix looks like
\begin{equation}
\begin{array}{lll}
r_{{\cal JB}}^{{\cal P}}\left( \left\{ \psi _{l},\varsigma _{1}\right\}
\right) & = & \psi _{1}\varsigma _{1}\left( A_{13}-\frac{1}{\varsigma _{1}}
A_{14}\right) \wedge B_{1}+\psi _{2}\varsigma _{1}H_{2,3}^{{\cal P}}\wedge
\left( A_{23}-\frac{1}{\varsigma _{1}}A_{24}\right) \\[0.3cm]
&  & \qquad +\; \psi _{1}\left( H_{1,4}^{{\cal P}}\wedge A_{1,4}+E_{1,2}\wedge
A_{2,4}\right) .
\end{array}
\end{equation}
The corresponding sumplectic form (in terms of
the algebra ${\mathfrak g}_{{\cal JB}}^{{\cal P}}$)
\begin{eqnarray*}
\omega _{{\cal JB}}^{{\cal P}}\left( \left\{ \psi _{i},\varsigma
_{1}\right\} \right) &=&-\sum_{l=1}^{2}\sum_{m=3}^{5-l}
\frac{1}{\alpha _{lm}}
E_{l,m}^{\ast }\left( \left[ ,\right] \right)  \\[0.3cm]
&=&-\left( \frac{1}{\psi _{1}\varsigma _{1}}E_{1,3}^{\ast }+\frac{1}{\psi
_{1}}E_{1,4}^{\ast }+\frac{1}{\psi _{2}\varsigma _{1}}
E_{2,3}^{\ast }\right)
\left( \left[ ,\right] \right) ,
\end{eqnarray*}
defines ${\mathfrak g}_{{\cal JB}}^{{\cal P}}$ as a Frobenius algebra.

The expression (\ref{tw-pjb}) for the enlarged peripheric twist
shows that the Hopf algebra $U_{{\cal B}}^{{\cal P}}(sl(4))$ can be
twisted further by the JF $e^{H_{1}^{\perp }\otimes \sigma _{1}}$:
\[
U(sl(4))\stackrel{{\cal F}_{{\cal B}_{1\prec 0}}^{{\cal P}}}
{\longrightarrow}
U_{{\cal B}}^{{\cal P}}(sl(4))\stackrel{e^{H_{1}^{\perp }\otimes
\sigma_{1}}}{\longrightarrow }U_{{\cal JB}}^{{\cal P}}(sl(4)).
\]
The final co-structure $\Delta _{{\cal JB}}^{{\cal P}}$ can be obtained
as $\Delta _{{\cal JB}}^{{\cal P}}=e^{{\rm ad}
\left( H_{1}^{\perp }\otimes \sigma_{1}\right) }
\circ \Delta _{{\cal B}}^{{\cal P}}$
and is defined by the set of coproducts:
\[
\begin{array}{lll}
\Delta _{{\cal JB}}^{{\cal P}}\left( E_{1,2}\right)  & = & E_{1,2}
\otimes e^{\sigma _{1,3}-\sigma _{1,4}}+1\otimes
E_{1,2}-H_{2,3}^{{\cal P}}\otimes E_{13}
e^{-\sigma _{2,3}}; \\[0.25cm]
\Delta _{{\cal JB}}^{{\cal P}}\left( E_{1,3}\right)  & = & E_{1,3}
\otimes e^{\sigma _{1,3}}+1\otimes E_{1,3}; \\[0.25cm]
\Delta _{{\cal JB}}^{{\cal P}}\left( E_{1,4}\right)  & = & E_{1,4}
\otimes e^{\sigma _{1,4}}+1\otimes E_{1,4}; \\[0.25cm]
\Delta _{{\cal JB}}^{{\cal P}}\left( E_{2,3}\right)  & = & E_{2,3}
\otimes e^{\sigma _{2,3}}+1\otimes E_{2,3}; \\[0.25cm]
\Delta _{{\cal JB}}^{{\cal P}}\left( E_{2,4}\right)  & = & E_{2,4}
\otimes e^{\sigma _{2,3}-\sigma _{1,3}}+e^{\sigma _{1,4}}
\otimes E_{2,4};  \\[0.25cm]
\Delta _{{\cal JB}}^{{\cal P}}\left( E_{3,4}\right)  & = &
E_{3,4}\otimes e^{-\sigma _{1,3}}+e^{\sigma _{1,4}}\otimes E_{3,4}
\\[0.2cm]
& & \qquad + H_{1}^{\perp }e^{\sigma _{1,4}}\otimes
 \left( e^{\sigma _{1,4}}-1 \right)
e^{-\sigma _{1,3}}+H_{2,3}^{{\cal P}}\otimes E_{24}e^{-\sigma _{2,3}};
\\[0.25cm]
\Delta _{{\cal JB}}^{{\cal P}}\left( H_{1,4}^{{\cal P}}\right)  & = &
H_{1,4}^{{\cal P}}\otimes e^{-\sigma _{1,4}}+1\otimes H_{1,4}^{{\cal P}
}-H_{1}^{\perp }\otimes \left( e^{-\sigma _{1,3}}-1\right)  \\[0.2cm]
& &\qquad -\left( e^{\sigma _{1,3}}-1\right)
\otimes E_{3,4}e^{-\sigma _{1,4}+
\sigma_{1,3}}-H_{1}^{\perp }\left( e^{\sigma _{1,3}}-1\right)
\otimes \left( 1-e^{-\sigma _{1,4}}\right) \\[0.20cm]
& &\qquad -\left( E_{1,2}+H_{2,3}^{{\cal P}}E_{1,3}\right)
\otimes E_{2,4}
e^{-\sigma_{1,4}-\sigma _{2,3}+\sigma _{1,3}};\\[0.25cm]
\Delta _{{\cal JB}}^{{\cal P}}\left( H_{2,3}^{{\cal P}}\right)  & = &
H_{2,3}^{{\cal P}}\otimes e^{-\sigma _{2,3}}+1\otimes H_{2,3}^{{\cal P}};
\\[0.25cm]
\Delta _{{\cal JB}}^{{\cal P}}\left( H_{1}^{\perp }\right)  & = &
H_{1}^{\perp }\otimes e^{-\sigma _{1,3}}+1\otimes H_{1}^{\perp };
\end{array}
\]
All the Cartan generators of $sl(4)$ participate in the
deformation that leads to this Hopf algebra
$U_{{\cal JB}}^{{\cal P}}(sl(4))$.

In the ordinary injection of the Poincar\'e algebra in $sl(4)$
(in the corresponding real form) the Cartan generators are
identified with the
diagonalizable rotation (usually $l_{12}$), boost ($n_{03}$)
and the dilatation operator.
Notice that here the commutative subalgebra of generators
$\{ A_{ij}=E_{ij}\; |\; i=1,2;\; j=3,4\}$ is just the subalgebra
of  Poincar\'e translations.
All of them have the especially simple quasi-primitive coproducts.
The reason is that in the twist ${\cal F}_{{\cal JB}}^{{\cal P}}$
the number of Jordanian factors ${\cal F}_{\cal J}$
is maximal and each of them
is attached to one of the ``momenta".

\subsect{$U\left( sl(7)\right) $}

Now we shall briefly expose the case of dimension seven to show some
peculiarities of odd dimension $N=2n-1$.
We also want the algebra $U(sl(N))$ to be
sufficiently large to contain nontrivially the generators of the type
$ C_{lm}\left( \left\{ \varsigma _{i}\right\} \right) $.

In this case the full peripheric chain has three links
$\left( n=4\; ;\; i,j=1,2,3\; ;\; z=N-n=3 \right) $,
\begin{equation}  \label{per-ch-sl7}
\begin{array}{rcl}
{\cal F}_{{\cal B}_{2\prec 0}}^{{\cal P}} & =
& \prod\limits_{k=0}^{\stackrel{2}{
\longleftarrow }}{\cal F}_{{\cal B}_{k}}^{{\cal P}}=
\prod\limits_{k=0}^{\stackrel{2 }{\longleftarrow }}
\left( \prod\limits_{s=k+2}^{6-k}e^{E_{k+1,s}\otimes
E_{s,7-k}}\right) e^{H_{k+1,7-k}^{{\cal P}}\otimes
\sigma _{k+1,7-k}} \\[5mm]
& = & e^{E_{3,4}\otimes E_{4,5}}e^{H_{3,5}^{{\cal P}}
\otimes \sigma_{3,5}} \\[2mm]
&  & \times \; e^{E_{2,3}\otimes E_{3,6}}e^{E_{2,4}
\otimes E_{4,6}}e^{E_{2,5}\otimes
E_{5,6}}e^{H_{2,6}^{{\cal P}}\otimes \sigma _{2,6}}\\[2mm]
&  & \times \; e^{E_{1,2}\otimes E_{2,7}}e^{E_{1,3}
\otimes E_{3,7}}e^{E_{1,4}\otimes
E_{4,7}}e^{E_{1,5}\otimes E_{5,7}}e^{E_{1,6}\otimes E_{6,7}}
e^{H_{1,7}^{ {\cal P}}\otimes \sigma _{1,7}}.
\end{array}
\end{equation}
The Cartan
generators are fixed according to the prescription (\ref{per-cart}):
\[
\begin{array}{rcl}
H_{k+1,N-k}^{{\cal P}} & = & \displaystyle \frac{1}{N}\left(
NE_{k+1,k+1}-\sum_{s=1}^{N}E_{s,s}\right) ,\qquad k=0,1,2 .
\end{array}
\]
The dimension of the subalgebra ${\bf H}^{\perp }$ in ${\bf H}(sl(7))$
that remains primitive after the twisting performed
 by the peripheric chain 
${\cal F}_{{\cal B}_{2\prec 0}}^{{\cal P}}$
\[
{\cal F}_{{\cal B}_{2\prec 0}}^{{\cal P}}:U(sl(7))
\longrightarrow U_{_{{\cal  B}}^{{\cal P}}}(sl(7))
\]
coincides with the dimension of ${\bf H}^{{\cal P}}$. We choose the
following three basic elements (see (\ref{add-gen})):
\begin{equation}
H_{i}^{\perp }=\frac{N-2i}{N}\sum_{l=1}^{N}E_{l,l}-
\sum_{m=i+1}^{N-i}E_{m,m},\qquad i=1,2,3 .
\end{equation}
It is easy to check that the subalgebras ${\bf H}^{{\cal P}}$ and
${\bf H}^{\perp}$ are primitive in the
twisted algebra $U_{{\cal B}}^{{\cal P}}(sl(N))$ :
\[
\begin{array}{lll}
\Delta _{{\cal B}_{2\prec 0}}^{{\cal P}}\left( H_{i}^{\perp }\right)
& = & H_{i}^{\perp }\otimes 1+1\otimes H_{i}^{\perp }, \\[0.25cm]
\Delta _{{\cal B}_{2\prec 0}}^{{\cal P}}\left( E_{j}^{{\cal P}}\right)
& = & E_{j}^{{\cal P}}\otimes 1+1\otimes E_{j}^{{\cal P}} ,
\end{array}
\qquad i,j=1,2,3,
\]
where $E_{j}^{{\cal P}}=E_{j,N-j}$. The generators
$\left\{ H_{i}^{\perp},E_{j}^{{\cal P}}\right\} $
form three (mutually commuting) Borel subalgebras:
\[
\left[ H_{i}^{\perp },E_{j}^{{\cal P}}\right] =
\delta _{ij}E_{j}^{{\cal P}}.
\]
This indicates that the peripheric chain
${\cal F}_{{\cal B}_{2\prec 0}}^{ {\cal P}}$ (\ref{per-ch-sl7})
can be enlarged by the factors
$\mit\Phi _{J_{i}}=e^{H_{i}^{\perp }\otimes \sigma _{i}}$,
with $\sigma _{i}=\ln \left(1+E_{i}^{{\cal P}}\right) $:
\[
{\cal F}_{{\cal JB}}^{{\cal P}}=\left( \prod\limits_{i=1}^{3}
e^{H_{i}^{\perp}\otimes \sigma _{i}}\right)
{\cal F}_{{\cal B}_{1\prec 0}}^{{\cal P}}.
\]
The corresponding universal ${\cal R}$--matrix has the form:
\[
{\cal R}_{{\cal JB}}^{{\cal P}}=\left( \prod\limits_{i=1}^{3}e^{\sigma
_{i}\otimes H_{i}^{\perp }}\right)
\left( {\cal F}_{{\cal B}_{2\prec 0}}^{{\cal P}}\right) _{21}
\left( {\cal F}_{{\cal B}_{2\prec 0}}^{{\cal P}
}\right) ^{-1}\left( \prod\limits_{i=1}^{3}e^{-H_{i}^{\perp }
\otimes \sigma_{i}}\right) .
\]
Notice that the additional JF's commute.

We have six deformation parameters: $\left\{ \psi _{l}\Rightarrow \xi
\psi _{l},\; \varsigma _{l}\;|\;l=1,2,3\right\} $. The expansion of
\\${\cal R}_{{\cal JB}}^{{\cal P}}
\left( \psi _{l},\varsigma _{l};\xi \right) $
with respect to $\xi $ exhibits the following classical $r$--matrix
\begin{equation}
\label{en-r-sl7}
\begin{array}{lll}
r_{{\cal JB}}^{{\cal P}}\left( \left\{ \psi _{l},\varsigma _{l}\right\}
\right)  & = & \psi _{1}\left[ H_{1}^{{\cal P}}\wedge
E_{1,7}+\sum_{k=2}^{5}E_{1,k}\wedge E_{k,7}+E_{1,6}\wedge \left(
E_{6,7}-\varsigma _{1}H_{1}^{\perp }\right) \right]   \\[0.2cm]
&  &\quad  +\; \psi _{2}\;\varsigma _{1}\left[ H_{2}^{{\cal P}}\wedge
E_{2,6}+\sum_{k=3}^{4}E_{2,k}\wedge E_{k,6}+E_{2,5}\wedge \left( E_{5,6}-
\frac{\varsigma _{2}}{\varsigma _{1}}H_{2}^{\perp }\right) \right]
\\[0.2cm]
&  & \quad +\; \psi _{3}\;\varsigma _{2}\left[ H_{3}^{{\cal P}}\wedge
E_{3,5}+E_{3,4}\wedge \left( E_{4,5}-\frac{\varsigma _{3}}{\varsigma _{2}}
H_{3}^{\perp }\right) \right] .
\end{array}
\end{equation}
The carrier algebra ${\mathfrak g}_{{\cal JB}}^{{\cal P}}$ of
the $r$--matrix $r_{{\cal JB}}^{{\cal P}}\left( \left\{ \psi _{l},
\varsigma _{l}\right\} \right) $ is
24D with the generators $\left\{ H_{i}^{{\cal P}},\right.$
$ H_{i}^{\perp},$ $\left.\left\{ E_{p,t}\;|\;p=1,2,3;\; t=2,\ldots ,7;
\;p<t\right\},\,
\left\{ E_{4,6},\;E_{4,7},\;E_{5,7}\right\} \right\} $.
This is an algebra of  motion
over the 12D space with translations
${\mathfrak g}_{{\cal P}  }=
\left\{ E_{p,t}\;|\;p=1,2,3;\; t=4,5,6,7\right\} $
and  the subalgebra
${\mathfrak g}_{{\cal H }}$ that is a direct sum of two 6D algebras
${\mathfrak g}_{{\cal H}}={\mathfrak g}_{ {\cal H}}^{\prime }\oplus
{\mathfrak g}_{{\cal H}}^{\prime\prime }$.
Each direct summand contain three Cartan
and three positive root generators:
${\mathfrak g}_{{\cal H} }^{\prime }=\left\{ H_{i}^{{\cal P}},
E_{p,t}\;|\;p,t=1,2,3;\right.$ $\left.p<t\right\} $ and
${\mathfrak g}_{{\cal H}}^{\prime \prime }=
\left\{ H_{i}^{\perp},\;E_{4,6}\,;E_{4,7},\;E_{5,7}\right\} $.
In these terms
${\mathfrak g}_{{\cal JB}}^{{\cal P}}$ is a semi-direct sum,
\begin{equation}
\label{g-struct}
{\mathfrak g}_{{\cal JB}}^{{\cal P}}=
\left({\mathfrak g}_{{\cal H}}^{\prime }
\oplus {\mathfrak g}_{{\cal H}}^{\prime \prime }\right) \vdash
{\mathfrak g}_{{\cal P}}.
\end{equation}
The image
$\varphi \left({\mathfrak g}_{{\cal JB}}^{{\cal P}}\right) $
of the map $ \varphi $ contains the generators
\[
\begin{array}{lllll}
B_{1}\left( \left\{ \varsigma _{1}\right\} \right)  & = & \displaystyle
\frac{1}{\varsigma
_{1}}E_{6,7}-H_{1}^{\perp }, &  & E_{i,j},\quad i,j=1,2,3, \\[0.3cm]
B_{2}\left( \left\{ \varsigma _{1,2}\right\} \right)  & = &
\displaystyle\frac{1}{\varsigma _{2}}E_{5,6}
+\frac{1}{\varsigma _{1}\varsigma _{2}}
E_{5,7}-H_{2}^{\perp }, &  & H_{i}^{{\cal P}}, \\[0.3cm]
B_{2}\left( \left\{ \varsigma _{1,2}\right\} \right)  & =
 &\displaystyle \frac{
1}{\varsigma _{3}}E_{4,5}
+\frac{1}{\varsigma _{2}\varsigma _{3}}E_{4,6}
+\frac{1}{\varsigma _{1}
\varsigma _{2}\varsigma _{3}}E_{4,7}
-H_{3}^{\perp }, &  & E_{i,7}. \\[0.3cm]
A_{i,4} & = &\displaystyle E_{i,4}+\frac{1}{\varsigma _{3}}E_{i,5}
+\frac{1}{\varsigma_{2}\varsigma _{3}}E_{i,6}
+\frac{1}{\varsigma _{1}\varsigma _{2}\varsigma
_{3}}E_{i,7}, &  & C_{4,6}  = \displaystyle E_{4,6}
+\frac{1}{\varsigma _{1}}E_{4,7}, \\[0.3cm]
A_{i,5} & = & \displaystyle E_{i,5}+\frac{1}{\varsigma _{2}}E_{i,6}
+\frac{1}{\varsigma_{1}\varsigma _{2}}E_{i,7}, &  & C_{4,7}  =  E_{4,6},
\\[0.3cm] A_{i,6} & = & \displaystyle E_{i,6}
+\frac{1}{\varsigma _{1}}E_{i,7},
&  & C_{5,7}  =\displaystyle  E_{5,7}-\frac{\varsigma _{2}}
{\varsigma_{3}}E_{4,7},
\end{array}
\]
In these terms the $r$--matrix $r_{{\cal JB}}^{{\cal P}}
\left( \left\{ \psi_{l},\varsigma _{l}\right\} \right) $
has the form
\begin{equation}
\label{r-pjb-7}
\begin{array}{rll}
&  & r_{\cal JB}^{\cal P}\left( \left\{ \psi _{l},\varsigma _{l}
\right\}\right) = \\[0.3cm]
 & &\qquad  \displaystyle\psi_{1}
\left[ H_{1}^{\cal P}\wedge A_{1,7}+E_{1,2}\wedge A_{2,7}
+E_{1,3}\wedge A_{3,7}+\left( A_{14}
-\frac{\varsigma_{2}}{\varsigma _{3}}
A_{15}\right) \wedge A_{47}\right. \\[0.3cm]
&  &\qquad +\left.\displaystyle \left( A_{15}
-\frac{\varsigma _{1}}{\varsigma _{2}} A_{16}\right)
\wedge \left( A_{57}+\frac{\varsigma _{2}}{\varsigma _{3}}
A_{47}\right)+\left( A_{16}-
\frac{1}{\varsigma _{1}} A_{17}\right) \wedge
\varsigma_{1}B_{1}\right]  \\[0.3cm]
&  &\qquad +\; \displaystyle\varsigma _{1}\psi _{2}
\left[ H_{2}^{\cal P}\wedge \left( A_{2,6}-
\frac{1}{\varsigma _{1}} A_{27}\right) +E_{2,3}\wedge
\left( A_{36}-\frac{1}{\varsigma _{1}} A_{37}\right) \right.
\\[0.3cm]
&  &\qquad +\;\displaystyle \left( A_{24}-\frac{\varsigma _{2}}
{\varsigma _{3}} A_{25}\right) \wedge
\left( A_{46}-\frac{1}{\varsigma _{1}} A_{47}\right) \\[0.3cm]
&  &\qquad + \left. \displaystyle \left( A_{25}
-\frac{\varsigma _{1}}{\varsigma _{2}} A_{26}\right)
\wedge \left( \frac{\varsigma _{2}}{\varsigma _{1}} B_{2}
-\frac{1}{\varsigma_{1}}\left( A_{57}
+\frac{\varsigma _{2}}{\varsigma _{3}}A_{47}\right)
\right) \right]  \\[0.3cm]
&  &\qquad \displaystyle
+ \; \varsigma _{2}\psi _{3}\left[ H_{3}^{\cal P}
\wedge \left( A_{3,5}-
\frac{\varsigma _{1}}{\varsigma _{2}}A_{36}\right)
+\left( A_{3,4}-\frac{\varsigma _{2}}{\varsigma _{3}}A_{35}\right)
\wedge \left( \frac{\varsigma _{3}}{\varsigma _{2}}B_{3}
-\frac{\varsigma _{1}}{\varsigma_{2}} A_{46}\right) \right] .
\end{array}
\end{equation}
This expression enables us to calculate the corresponding form
$\omega$ that has a
simple structure in terms of the algebra
${\mathfrak g}_{{\cal JB}}^{{\cal P}}$:
\[
\omega _{{\cal JB}}^{{\cal P}}\left( \left\{ \psi _{l},\varsigma
_{l}\right\} \right) =-\sum_{l=1}^{3}
\sum_{m=4}^{8-l}\frac{1}{\alpha _{lm}}E_{l,m}^{\ast }
\left( \left[ ,\right] \right) ;\quad \alpha _{lm}=\psi
_{l}\varsigma _{\left( N-m\right) }.
\]

Thus, when the Hopf algebra
$U_{_{{\cal B}_{2\prec 0}}^{{\cal P}}}(sl(7))$ is
twisted by the additional JF's
$\prod\limits_{i=1}^{3}e^{H_{i}^{\perp }\otimes \sigma _{i}}$, the result
\[
\prod\limits_{i=1}^{3}e^{H_{i}^{\perp }\otimes \sigma _{i}}:U_{_{{\cal B}
_{2\prec 0}}^{{\cal P}}}(sl(7))\longrightarrow U_{_{{\cal JB}}^{{\cal P}
}}(sl(7))
\]
is a quantization of the initial Hopf algebra $U(sl(7))$ in the direction
defined by the $r$--matrix $r_{{\cal JB}}^{{\cal P}}
\left( \left\{ \psi_{l},\varsigma _{l}\right\} \right) $
((\ref{en-r-sl7}) or (\ref{r-pjb-7})) with  carrier algebra
${\mathfrak g}_{{\cal JB}
}^{{\cal P}}$ (see (\ref{g-struct})).

\subsect{$U\left( sl(3)\right)$}

To conclude the examples it is worth mentioning the degenerate case of
the enlarged peripheric chain, the PET ${\cal F}_{{\cal E}}^{{\cal P}
}=e^{E_{1,2}\otimes E_{2,3}}e^{H^{{\cal P}}\otimes \sigma _{1,3}}$ in $
U(sl(3))$ equipped with the additional Jordanian factor:
\begin{equation}
\label{en-per-sl3}
{\cal F}_{{\cal JE}}^{{\cal P}}=e^{H_{13}^{\perp }\otimes
\sigma _{12}\left(\varsigma \right) }
e^{\psi E_{1,2}\otimes E_{2,3}}e^{H^{{\cal P}}\otimes
\sigma _{1,3}\left( \psi \right) },
\end{equation}
(see (\ref{simp-add-jord})).
Here the generators $H^{{\cal P}}$ and $H_{13}^{\perp }=H^{\perp }$ were
defined in section~\ref{deformationcarriers}, formula (\ref{generat-3}).
The corresponding classical $r$--matrix (\ref{sl3-pr})
is evidently one-parametric:
\begin{equation}
\label{r-en-per-sl3}
\begin{array}{lll}
r_{{\cal JE}}^{{\cal P}}\left( \varsigma \right)  & = & H^{{\cal P}}\wedge
E_{13}+E_{12}\wedge \left( E_{23}-\varsigma H_{13}^{\perp }\right) .
\end{array}
\end{equation}
Its carrier algebra ${\mathfrak g}_{{\cal JE}}^{{\cal P}}$ is 4D with two
Cartan and two ``translation'' generators: $\left\{ H^{{\cal P}} \right.
,\; H_{13}^{\perp },\; E_{12},$
$ \left. E_{13}\right\} $. The natural basis in the
image $\varphi \left( {\mathfrak g}_{ {\cal JE}}^{{\cal P}}\right) $
of the map $\varphi $ is formed by the elements:
\[
\begin{array}{lllll}
B\left( \left\{ \varsigma \right\} \right)  & = & \frac{1}{\varsigma }
E_{23}-H_{13}^{\perp }, &  & H^{{\cal P}}, \\[0.3cm]
A_{12} & = & E_{12}+\frac{1}{\varsigma }E_{13}, &  & E_{13}.
\end{array}
\]
The $\omega $--form corresponding to the $r$--matrix (\ref{r-en-per-sl3}) is
\[
\omega _{{\cal JE}}^{{\cal P}}\left( \varsigma \right)
 =-\frac{1}{\varsigma }  E_{12}^{\ast }
\left( \left[ ,\right] \right) -E_{13}^{\ast }
\left( \left[ , \right] \right) .
\]

The important fact here is that for the 4D carrier
algebra ${\mathfrak g}_{{\cal JE}}^{{\cal P}}$
we can easily compose the other twisting element that
contains no external generators. Consider the basic generators
$\left\{ H_{12}^{{\cal \perp }},
H_{13}^{\perp },E_{12},E_{13}\right\} $ correlated with
the direct sum structure
${\mathfrak g}_{{\cal JE}}^{{\cal P}}={\bf B}_{2}
\left( H_{12}^{  {\cal \perp }},E_{13}\right)
\oplus {\bf B}_{2}\left( H_{13}^{{\cal \perp } },
E_{12}\right) $.
This is obviously the two-Jordanian situation and the
corresponding twist is
\[
{\cal F}_{{\cal JJ}}=e^{H_{13}^{\perp }\otimes
\sigma _{12}\left( \varsigma \right) }
e^{H_{12}^{{\cal \perp }}\otimes \sigma _{13}\left( \psi \right) }.
\]
Notice that the $r$--matrix
\[
r_{{\cal JJ}}\left( \varsigma \right) =H_{12}^{{\cal \perp }}\wedge
E_{13}+\eta H_{13}^{\perp }\wedge E_{12}
\]
is equivalent to (\ref{r-en-per-sl3}).
In this particular case we can construct
two different twists,
${\cal F}_{{\cal JJ}}$ and ${\cal F}_{{\cal JE}}^{\cal P}$,
that quantize the $r$--matrix
$r_{{\cal JJ}} \approx r_{{\cal JE}}^{\cal P}$.
One of them has the carrier algebra
${\mathfrak g}_{{\cal JE}}^{{\cal P}}$, the
other depends on the full set of generators of ${\cal B}^+(sl(3))$.
In other words, the twist ${\cal F}_{{\cal JE}}^{\cal P}$
uses not only the adjoint
operators of the elements of
${\mathfrak g}_{{\cal JE}}^{{\cal P}}$
but also the morphisms
external with respect to  ${\mathfrak g}_{{\cal JE}}^{{\cal P}}$.

This is the explicit demonstration of the fact that
there exist different quantizations
of one and the same Lie-Poisson structure.
\sect{Conclusions}

We have explicitly quantized the sets of Lie-Poisson structures
defined on the groups of motion. When they are of the
$r_{{\cal JB}}^{{\cal P}}$--type the
carrier algebras of the $r$-matrix and the twist
${\cal F}_{{\cal JB}}^{{\cal P}} \left(
\left\{ \nu _{l},\varkappa _{l},\rho_{i}\right\} \right) $
are different (${\mathfrak g}_{{\cal JB}}^{{\cal P}}$
and ${\bf B}^{+}\left( sl(N)\right) $, respectively).
It is necessary to stress that the defining space of the twisting element
${\cal F}_{{\cal JB}}^{{\cal P}}\left(
\left\{ \nu _{l},\varkappa _{l},\rho_{i}\right\} \right) $
cannot be restricted to the subspace
$U\left( {\mathfrak  g}_{{\cal JB}}^{{\cal P}}\right)
\otimes U\left( {\mathfrak g}_{{\cal JB}}^{{\cal P}}\right) $.
The quantization of the
$r$--matrix $r_{{\cal JB}}^{{\cal P}}$ with  carrier subalgebra
${\mathfrak g}_{{\cal JB}}^{{\cal P}}$
was performed by the twist ${\cal F}_{{\cal JB}}^{{\cal P}
}\left( \left\{ \nu _{l},\varkappa _{l},\rho _{i}\right\} \right) $ with
carrier ${\bf B}^{+}\left( sl(N)\right) $.
On the other hand, it is  well
known \cite{D83} that to find the twisting element
it is sufficient to know
the classical $r$--matrix, there exists the solution
of the Drinfeld equations
corresponding to $r$ and defined over
${\mathfrak g}_{{\cal JB}}^{{\cal P}}$.
It was demonstrated above that this is not the only way to perform
the explicit twisting. The alternative possibility
is to inject the carrier into the larger
algebra and to search for the solutions there. Thus, to perform
the quantization corresponding to the carrier
${\mathfrak g}_{{\cal JB}}^{{\cal  P}}$ we
have used the enlarged peripheric twist
${\cal F}_{{\cal JB}}^{{\cal P}}\left(
\left\{ \nu _{l},\varkappa _{l},\rho _{i}\right\} \right) $
that includes the adjoint operators external with respect to
${\mathfrak g}_{{\cal JB}}^{{\cal P}}$.
This method can be used to facilitate the solution of the
Drinfeld equations: in the
nonordinary cases it might
be reasonable to enlarge the algebra and thus
include the external morphisms.

In the general situation the explicit form of the twisting element
with the carrier ${\mathfrak g}_{{\cal JB}}^{{\cal P}}$
(that depends only on the generators of
${\mathfrak g}_{{\cal JB}}^{{\cal P}}$ )
is unknown to us.  The only exception is  the
degenerate case of $U\left( sl(3)\right) $ considered in the
end of the previous Section. The same properties can be
found when $N$ is arbitrary but
the chain is degenerate and consists of only one link.
In this single link the
extension can contain $N-2$  basic extending
factors and (when the whole link is peripheric) can
lead to $N-2$ additional primitive generators.
The enlarged peripheric twist
will have $N-2$ additional Jordanian factors.
At the same time here the
carrier ${\mathfrak g}_{{\cal JE}}^{{\cal P}}$
is obviously equivalent to the direct sum
of $N-1$ commuting Borel subalgebras and the corresponding
multi-Jordanian twist can be composed.

The algebras ${\mathfrak g}_{{\cal JB}}^{{\cal P}}=
{\mathfrak g}_{{\cal H}}\vdash {\mathfrak g}_{{\cal P}}$ are   algebras of
motion of   special type. They have the $(z\left( N-z\right))$-dimensional
subspace corresponding to the translations
$A_{lm}\in {\mathfrak g}_{{\cal P}}$
. The subalgebra ${\mathfrak g}_{{\cal H}}$ contains the Borel subalgebra
${\bf B}^{+}\left( gl\left( z\right) \right) $, which is the multidimensional
analogue of the dilatation operator in the conformal algebra,
and the subalgebra ${\bf C}^{+}\left( gl\left( n-1\right) \right) $
that can be treated as a contraction of
${\bf B}^{+}\left( gl\left( n-1\right) \right) $. The Cartan
subalgebra in ${\mathfrak g}_{{\cal H}}$ contains the ordinary set of
simultaneously diagonalizable operators in the algebra of motion
${\mathfrak g}_{{\cal JB}}^{{\cal P}}$ including
the subset of dilatation-like operators $B_i$. The
applications of the enlarged peripheric chains of twists
studied in this paper to quantize Lie-Poisson structures defined on groups
of motion could be of considerable interest.

To study the peripheric chains for simple Lie algebras of series $A_n$ we have
used the auxiliary Reshetikhin ``rotation" in the root space. This method is
quite general and can be applied to an arbitrary simple algebra.

\section*{Acknowledgments}

One of the authors (V. L.) would like to thank the Vicerrectorado de
Investigaci\'on de la Universidad de Valladolid for supporting his
stay. This work has also  been partially supported by  DGES of the
Ministerio de Educaci\'on y Cultura de Espa\~na under Project PB98--0360, the
Junta de  Castilla y Le\'on (Spain) and by the
Russian Foundation for Fundamental Research under the grant 00--01--00500.



\end{document}